\documentclass[10pt]{article}
\usepackage{amssymb,amsmath,textcomp}
\newtheorem{satz}{Theorem}[section]

\newtheorem{lemma}{Lemma}[section]

\newtheorem{bemerk1}{Remark}[section]

\newtheorem{bei}{Example}[section]
\newtheorem{korollar}{Corollary}[section]
\newtheorem{prop}{Proposition}[section]
\newcommand{\iR}{\mathbb{R}}
\newcommand{\iN}{\mathbb{N}}

\newcommand{\iC}{\mathbb{C}}

\newcommand{\oH}{\hspace*{0.39em}\raisebox{0.6ex}{\textdegree}\hspace{-0.72em}H}

\begin{document}
\begin{center}
{\bf\Large Decay estimates for time-fractional and other non-local in time subdiffusion equations
in $\iR^d$}
\end{center}
\vspace{0.7em}
\begin{center}
Jukka Kemppainen\footnote{J.K., J.S. and R.Z. were partially supported by Academy of Finland project 138738 (Nonlinear PDEs, project lead by professor Juha Kinnunen)}, Juhana Siljander\footnote{J.S. was supported by Academy of Finland postdoctoral research grant 259363}, Vicente Vergara\footnote{V.V.\ was partially supported by FONDECYT grant 1110033.}, and Rico Zacher\footnote{R.Z.\ was supported by a Heisenberg fellowship of the German Research Foundation (DFG), GZ  Za 547/3-1.}
\end{center}
\begin{abstract}
We prove optimal estimates for the decay in time of solutions to a rather general class of non-local in time subdiffusion
equations in $\iR^d$. An important special case is the time-fractional diffusion equation, which has seen much interest
during the last years, mostly due to its applications in the modeling of anomalous diffusion processes. We follow three different approaches and techniques to study this particular case: (A) estimates based on the fundamental solution and Young's inequality, (B)
Fourier multiplier methods, and (C) the energy method. It turns out that the decay behaviour is markedly different from the
heat equation case, in particular there occurs a {\em critical dimension phenomenon}. The general subdiffusion case is treated by method (B) and relies on a careful estimation of the underlying relaxation function. Several examples of kernels,
including the ultraslow diffusion case, illustrate our results.
\end{abstract}
\vspace{0.7em}
\begin{center}
\end{center}

\section{Introduction and main results}
We study the temporal decay of solutions to the non-local in time diffusion
 equation
\begin{equation} \label{subdiffintro}
\partial_t \big(k\ast [u-u_0]\big)-\Delta u=0,\quad
t>0,\,x\in \iR^d,
\end{equation}
together with initial condition
\begin{equation} \label{subdiffinitintro}
u|_{t=0}=u_0,\,x\in \iR^d.
\end{equation}
Here $u_0(x)$ is a given datum. $k\ast v$ denotes the convolution on the positive halfline $\iR_+:=[0,\infty)$ w.r.t.\ the time variable,
that is
$(k\ast v)(t)=\int_0^t k(t-\tau)v(\tau)\,d\tau$, $t\ge 0$.
We assume that $k$ is a kernel of type $\mathcal{PC}$, by which we mean that the following condition is satisfied.
\begin{itemize}
\item [{\bf ($\mathcal{PC}$)}] $k\in L_{1,\,loc}(\iR_+)$ is nonnegative and nonincreasing, and there exists a kernel $l\in L_{1,\,loc}(\iR_+)$ such that
$k\ast l=1$ in $(0,\infty)$.
\end{itemize}
In this case we also write $(k,l)\in \mathcal{PC}$. Note that $(k,l)\in {\cal PC}$ implies that $l$ is completely
positive, cf.\ \cite[Theorem 2.2]{CN}, in particular $l$ is nonnegative.

An important example for a pair $(k,l)\in \mathcal{PC}$ is given by $(k,l)=(g_{1-\alpha},g_\alpha)$ with $\alpha\in(0,1)$, where $g_\beta$ denotes the standard kernel
\[
g_\beta(t)=\,\frac{t^{\beta-1}}{\Gamma(\beta)}\,,\quad
t>0,\quad\beta>0.
\]
In this case, the term $\partial_t(k\ast v)$ becomes the Riemann-Liouville fractional derivative
$\partial_t^\alpha v$ of order $\alpha$ (cf.\ \cite{KST}) and (\ref{subdiffintro}) is called {\em time-fractional diffusion equation}.

Another interesting example is given by the pair
\begin{equation} \label{ultrapair}
k(t)=\int_0^1 g_\beta(t)\,d\beta,\quad l(t)=\int_0^\infty \,\frac{e^{-st}}{1+s}\,{ds},\quad t>0.
\end{equation}
In this situation the operator $\partial_t(k\ast \cdot)$ is a so-called operator of
{\em distributed order} and (\ref{subdiffintro}) is an example of a so-called {\em ultraslow diffusion equation}, see e.g.\ \cite{Koch08}.

\smallbreak

{\bf Applications.} Problems of the form (\ref{subdiffintro}), in particular the time-fractional diffusion equation,
have seen much interest during the last years, mostly due to their applications in the modeling of anomalous diffusion, see
e.g.\ \cite{Koch08}, \cite{Koch11}, \cite{Metz}, \cite{Uch} and the references therein for the physical background.
To provide some more specific motivation, let $Z(t,x)$ denote the fundamental solution of (\ref{subdiffintro}) satisfying  $Z|_{t=0}=\delta_0$. If $k$ is
a kernel of type ($\mathcal{PC}$), then this fundamental solution can be constructed via subordination from the heat kernel  and one can show that $Z(t,\cdot)$ is a probability density function (pdf) on $\iR^d$ for all $t>0$, see the proof of \cite[Theorem 3]{Koch11} and Section \ref{SecFS} below. Given a pdf $u_0$ on $\iR^d$ satisfying some appropriate conditions, the solution $u(t,\cdot)$ of the initial-value problem (\ref{subdiffintro}), (\ref{subdiffinitintro}) is given by
the formula (\ref{solutionformula}) below, and thus is a pdf on $\iR^d$  for all $t>0$. So in this case, the problem (\ref{subdiffintro}), (\ref{subdiffinitintro}) describes the evolution of a pdf on $\iR^d$.

An important quantity that measures the dispersion of random processes and that
describes how fast particles diffuse is the {\em mean
square displacement}. It can be determined in experiments and is defined (in our situation) as
\begin{equation} \label{MSDdef}
m(t)=\int_{\iR^d}|x|^2 Z(t,x)\,dx,\quad t>0.
\end{equation}
In the case of the classical diffusion equation (i.e.\ $\alpha=1$) $m(t)=ct$, $t>0$,
with some constant $c>0$. In the time-fractional diffusion case (i.e.\ the first
example) one observes that $m(t)=ct^\alpha$ (cf.\ \cite{Metz}), which shows that the diffusion is slower than in the classical case of Brownian motion. During the recent decades, experimental studies have shown that
there is an abundance of processes that have such a {\em power-law} mean square displacement, see \cite{BG}, \cite{Metz}, \cite{Metz2}, \cite{Uch} and the references given therein. An important application is the diffusion on fractals
like e.g.\ some amorphous semiconductors \cite{Metz}, \cite{Uch}.
In our second example, the mean square displacement $m(t)$ behaves like $c\log t$ for $t\to \infty$, see \cite{Koch08}. In this case
(\ref{subdiffintro}) describes a so-called {\em ultraslow diffusion} process. Such processes have been extensively studied recently. They appear, for example, in polymer physics \cite{SchSokBl}, diffusion in disordered media ({\em Sinai diffusion}) \cite{Sinai}, and in diffusion generated by iterated maps \cite{DrKl}.
We would like to point out that in our setting, that is, with a pair of kernels $(k,l)\in \mathcal{PC}$
the mean square displacement is given by
\begin{equation} \label{msdformula}
m(t)= \,2d\,(1\ast l)(t),\quad t> 0,
\end{equation}
see Lemma \ref{MSD}.

Another context where equations of the form (\ref{subdiffintro}) and nonlinear variants of them appear is the modeling
of dynamic processes in materials with {\em memory}. Examples are given by the theory of heat
conduction with memory, see e.g.\ \cite{JanI} and the references therein as well as \cite{VZ}, and the
diffusion of fluids in porous media with memory, cf.\ \cite{CapuFlow}, \cite{JakuDiss}.

\medbreak

\noindent In view of condition ($\mathcal{PC}$) the problem (\ref{subdiffintro}), (\ref{subdiffinitintro}) can be reformulated as an abstract Volterra
equation on the positive halfline with a completely positive kernel; this can be seen by convolving the PDE with the kernel
$l$. There has been a substantial amount of work on such abstract Volterra and integro-differential equations since the 1970s,
in particular on existence and uniqueness, regularity, and long-time behaviour of solutions,
see, for instance, \cite{CNa}, \cite{CN}, \cite{Grip1}, \cite{ZEQ}, and the monograph \cite{JanI}.

One of the main objectives of this paper is to prove sharp estimates for the temporal decay of solutions to (\ref{subdiffintro}),
(\ref{subdiffinitintro}). We point out that for non-local in space diffusion equations, in particular space-fractional diffusion equations, corresponding results have been obtained recently, see e.g.\ \cite{BrdP},\cite{CaffVazq11}, \cite{CCR},
\cite{IR}, \cite{Vazq14}. Concerning non-local in time diffusion, the case of a bounded domain with homogeneous Dirichlet boundary condition and a kernel of type $\mathcal{PC}$ has been studied recently in \cite{VZ}. The decay estimates obtained in
\cite{VZ} are optimal and even true in the case of a uniformly elliptic operator
in divergence form with rough coefficients, cf.\ the remarks following formula (\ref{udecay}) below. For the time-fractional case with Laplacian we also refer to
\cite{MNV} and \cite{NSY}. Moreover, decay estimates for the time-fractional diffusion equation in $\iR^d$ have also gathered interest in the engineering community \cite{MZL}. However, due to the unexpected critical dimension phenomenon
(to be discussed below) some of the subtleties of the equation have been overlooked there.

Under appropriate conditions on $u_0$ the solution to (\ref{subdiffintro}),
(\ref{subdiffinitintro}) is given by
\begin{equation} \label{solutionformula}
u(t,x)=\int_{\iR^d} Z(t,x-y) u_0(y)\,dy,\quad t>0,\,x\in \iR^d.
\end{equation}
Given $u_0$ in some Lebesgue space $L_q(\iR^d)$, $q\in [1,\infty]$ we want to understand the decay behaviour
of $|Z(t,\cdot)\star u_0|_{L_r(\iR^d)}$ as $t\to \infty$ for suitable $r\in [1,\infty]$; here $f_1\star f_2$ denotes the convolution of $f_1,\,f_2$ on $\iR^d$. It turns out that the situation is markedly different from that in the case of the heat
equation, where $Z(t,x)=H(t,x)=(4\pi t)^{-d/2}\exp(-{|x|^2/4t})$ is the Gaussian heat kernel.

\smallbreak

{\bf Time-fractional diffusion.} Let us first consider for simplicity the time-fractional diffusion equation (i.e.\ $k=g_{1-\alpha}$, $\alpha\in (0,1)$) and the case $r=2$. Given $u_0\in L_2(\iR^d)$ we do not have in general any decay neither for $|Z(t,\cdot)\star u_0|_{L_2(\iR^d)}$ nor
for $|H(t,\cdot)\star u_0|_{L_2(\iR^d)}$. Now suppose that $u_0\in L_2(\iR^d)\cap L_1(\iR^d)$. Then it is well-known that $u(t,\cdot):=H(t,\cdot)\star u_0$ decays in the $L_2$-norm as
\[
|u(t,\cdot)|_2 \lesssim t^{-\frac{d}{4}},\quad t>0,
\]
and this estimate is the best one can obtain in general (see e.g.\ \cite{BS}). Here $|v|_2:=|v|_{L_2(\iR^d)}$, and $v(t)\lesssim w(t),\,t>0$ means that there exists a constant $C>0$ such that $v(t)\le Cw(t)$, $t>0$. In the case of time-fractional diffusion $u(t,\cdot):=Z(t,\cdot)\star u_0$ exhibits the decay behaviour
\begin{equation} \label{L2intro}
|u(t,\cdot)|_2 \lesssim t^{-\min\{\frac{\alpha d}{4},\alpha\}},\quad t>0,\;d\in \iN\setminus \{4\},
\end{equation}
see Corollary \ref{L2estimate} below.
Whereas for the heat equation the decay rate increases with the dimension $d$, time-fractional diffusion leads to the
phenomenon of a {\em critical dimension}, which is $d=4$ in this case. Below the critical dimension the rate increases with
$d$, the exponent being $\alpha$ times the one from the heat equation, while above the critical dimension the decay rate
is the same for all $d$, namely $t^{-\alpha}$. The reason why the decay rate does not increase any further with $d$ lies in the fact that $t^{-\alpha}$ (up to a constant) coincides with the decay rate in the case of a bounded domain and homogeneous Dirichlet boundary condition, cf.\ the remarks on problem (\ref{problembdddom}) below. This reveals another interesting phenomenon: In the time-fractional case the diffusion is so slow that in higher dimensions ($d$ above the critical dimension) restriction to a bounded domain and the requirement of a homogeneous Dirichlet boundary condition do not improve the rate of decay. This phenomenon cannot be
observed in the classical diffusion case, where we always have exponential (and thus a better) decay in the case of a bounded domain.

We point out that the estimate (\ref{L2intro}) is the best one can get in general, see Theorem \ref{theoremlower} and Example \ref{beifrac} below. In the case of the critical dimension
$d=4$ we obtain the same decay rate as for $d>4$, however we have to replace the $L_2$-norm by the weak $L_2$-norm,
that is we find that
\[
|u(t,\cdot)|_{2,\infty} \lesssim t^{-\alpha},\quad t>0.
\]

In the more general case where $r\in (1,\infty)$ and $u_0\in L_1(\iR^d)\cap L_r(\iR^d)$, the critical dimension (which is in general not an integer) is given by
\[
d_{crit}=\frac{2r}{r-1}.
\]
Assuming $d\ge 3$ we show that
$|u(t,\cdot)|_r  \lesssim\, t^{ -\frac{\alpha d }{2}\,\left(1-\frac{1}{r}\right)}$, if $d<d_{crit}$, and
$|u(t,\cdot)|_{r,\infty} \lesssim\,t^{-\alpha}$, if $d=d_{crit}$, as well as
$|u(t,\cdot)|_r \lesssim\,t^{-\alpha}$, whenever $d>d_{crit}$. Moreover, in the case $d<3$ we always have a subcritical decay behaviour, see Corollary \ref{Lrkorollar}.

The critical dimension phenomenon can be well understood by looking at the asymptotic properties of the fundamental
solution $Z$. In the time-fractional case, $Z$ can be represented by Fox-functions and corresponding estimates are known (cf.\ \cite{Koch} and Prop.\ \ref{Kochubei1}
below). They show that for $d\ge 2$ the fundamental solution $Z(t,x)$ does not only have a singularity at $t=0$ but also at
the origin $x=0$, which is a striking difference to the heat kernel. If, for example, $d\ge 3$ and $t^{-\alpha}|x|^2\le 1$,
we have the sharp estimate $Z(t,x)\le C t^{-\alpha} |x|^{-d+2}$, thus $Z$ is (up to a constant) bounded above by $g_{1-\alpha}(t)$ times
the Newtonian potential w.r.t.\ $x$. This provides some interesting insight into how the fundamental solutions from the
Poisson and the heat equation interpolate in case of fractional dynamics.

Relying on the asymptotic bounds for (the time-fractional) $Z$, in Section \ref{fractionalheat} we first derive estimates for $|Z(t,\cdot)|_p$ (for suitable $p\in [1,\infty]$), which by Young's inequality then yield bounds for $|Z(t,\cdot)\star u_0|_r$. We further look at gradient estimates for $Z(t,\cdot)$ and $Z(t,\cdot)\star u_0$, and we show that for integrable initial data $u_0$ the asymptotic behaviour of $Z(t,\cdot)\star u_0$ as $t\to \infty$ is described by a multiple of $Z(t,x)$. This is the time-fractional analogue of a well-known result for the heat equation, cf.\ \cite[Prop.\ 48.6]{QS}.

\smallbreak

{\bf The general subdiffusion case.} Turning to the more general subdiffusion equation (\ref{subdiffintro}) with $k$ being a kernel of type $\mathcal{PC}$, we first note that asymptotic bounds for $Z(t,x)$ like those in the time-fractional case do not seem to be known in the general
case. However, we mention \cite{Koch08} where the ultraslow-diffusion case was studied. To derive decay estimates
for $Z(t,\cdot)\star u_0$ we develop a theory which is based on tools from harmonic analysis and a careful estimation of the Fourier symbol $\tilde{Z}(t,\xi)$ of $Z$ w.r.t.\ the spatial variable.

Taking the Fourier transform w.r.t.\ $x$ (defined as in (\ref{FourierDef}) below) we see that
$\tilde{Z}(t,\xi)$ solves the problem
\[
\partial_t \big(k\ast [\tilde{Z}-1]\big)+|\xi|^2 \tilde{Z}=0,\quad
t>0,\,\xi\in \iR^d,\quad \tilde{Z}(0,\xi)=1,\quad \xi\in \iR^d,
\]
that is, we have
\[
\tilde{Z}(t,\xi)=s(t,|\xi|^2),\quad t\ge 0,\;\xi\in \iR^d.
\]
Here the so-called {\em relaxation function} $s(t,\mu):=s_\mu(t)$, $t\ge 0$, is defined for the parameter $\mu\ge 0$ as the solution of the Volterra integral equation
\begin{equation} \label{smudef}
s_\mu(t)+\mu(l\ast s_\mu)(t)=1,\quad t\ge 0.
\end{equation}
Note that $s_0\equiv 1$ and that (\ref{smudef}) is equivalent to the integro-differential equation
\[
\frac{d}{dt}\,\left(k\ast [s_\mu-1]\right)(t)+\mu s_\mu(t)=0,\quad t>0,\quad s_\mu(0)=1.
\]
It is known that the assumption $(k,l)\in \mathcal{PC}$ implies that $s_\mu$ is nonnegative, nonincreasing, and that $s_\mu\in H^1_{1,\,loc}(\iR_+)$; moreover
$\partial_\mu s_\mu(t)\le 0$, see e.g.\ Pr\"uss \cite{JanI}. Furthermore, it has been shown
in \cite{VZ} that for any $\mu\ge 0$ there holds
\begin{equation} \label{smubounds}
\frac{1}{1+\mu \,k(t)^{-1}}\,\le s_\mu(t)\le \,\frac{1}{1+\mu \,(1\ast l)(t)}\,,\quad
\mbox{a.a.}\;t>0,
\end{equation}
which also shows that
\[
\big[1-s_\mu(t)\big] k(t)\,\le \mu s_\mu(t)\le \,\big[1-s_\mu(t)\big]\,\frac{1}{(1\ast l)(t)} \,,\quad \mbox{a.a.}\;t>0.
\]
This implies that for any fixed $\mu>0$, $s_\mu(t)$ cannot decay faster than the kernel $k(t)$, and $s_\mu(t)$ decays
at least like $(1\ast l)(t)^{-1}$. Note that $\lim_{t\to \infty} s_\mu(t)=0$ if and only if $l\notin L_{1}(\iR_+)$, see e.g.\
\cite[Lemma 6.1]{VZ}.

By means of (\ref{smubounds}) and Plancherel's theorem we are able to prove the following $L_2$-decay estimate.
Suppose $u_0\in L_1(\iR^d)\cap L_2(\iR^d)$, then
\begin{equation} \label{L2I}
|Z(t,\cdot)\star u_0|_2 \lesssim \,\big[(1\ast l)(t)\big]^{-\min\{1,\frac{d}{4}\}},\quad t> 0,\;d\in \iN\setminus\{4\},
\end{equation}
see Theorem \ref{thmupper}. We also show that this estimate gives the optimal rate of decay provided $(1\ast l)(t)^{-1}\le
C k(t)$ for $t\ge T$, where $C,T$ are some positive constants. This applies in particular to the time-fractional diffusion equation, where $k(t)$ and $(1\ast l)(t)^{-1}=g_{1+\alpha}(t)^{-1}$ decay like $c t^{-\alpha}$, and to the
ultraslow diffusion equation, where
$k(t)$ and $(1\ast l)(t)^{-1}$ decay like $c (\log t)^{-1}$. (\ref{L2I}) also reveals that the phenomenon of critical dimension extends to the general case of a kernel $k$ of type $\mathcal{PC}$.

In order to obtain sharp $L_r$-decay estimates and to treat the critical dimension case we go a step further and
derive bounds for the partial derivatives of arbitrary order of the functions $\xi\to |\xi|^{2\delta}\tilde{Z}(t,\xi)$, where $t>0$ is fixed and $\delta\in (0,1]$ is a parameter. These bounds are uniform w.r.t.\ time and allow us to apply Mihlin's multiplier theorem.
This leads to the desired $L_r$- (and by interpolation also to the $L_{r,\infty}$-) decay estimates, which generalize those from the time-fractional case.
Here in our analysis we make use of another important property of the relaxation function $s(t,\mu)$ which is the {\em complete monotonicity w.r.t.\ to the parameter $\mu$} for all $t\ge 0$, see Section \ref{SecLR} below. This property has already been known,
e.g.\ it follows from results in \cite[Section 4]{JanI}. In the present paper we provide a new proof of this property, which is rather short and relies on a comparison principle for a certain type of integro-differential equation.

By switching the kernels from the ultraslow diffusion case one obtains an interesting example of a pair $(k,l)\in {\cal PC}$ where
$k(t)$ decays like $t^{-1}$ and $(1\ast l)(t)^{-1}$ like $c t^{-1}\log t$. Here the upper bound in (\ref{smubounds}) does
not lead to the optimal decay. However, it has been shown in \cite{VZ} that $s_\mu(t)\le \frac{C}{1+\mu t}$ for all $t,\mu\ge 0$, and thus a simple modification of our original proof yields the optimal $L_2$-estimate
\begin{align*}
|u(t,\cdot)|_2  & \lesssim t^{-\min\{1,\frac{d}{4}\}},\quad t> 0,\;d\in \iN\setminus\{4\},\\
|u(t,\cdot)|_{2,\infty} & \lesssim t^{-1},\quad t> 0,\;d=4,
\end{align*}
provided $u_0\in L_1(\iR^d)\cap L_2(\iR^d)$, cf.\ Example \ref{exotic}. This is a remarkable result as for $d\le 3$ the decay rate is the same as for the heat equation!

\smallbreak

{\bf Known results in the bounded domain case.} It is instructive to compare our decay results with what is known in the case of a bounded domain and a homogeneous
Dirichlet condition. Let $\Omega\subset \iR^d$ be a bounded domain, $u_0\in L_2(\Omega)$, and assume again that
$k$ is of type $\mathcal{PC}$. We consider
the problem
\begin{align}
\partial_t \big(k\ast [u-u_0]\big)-\Delta u & =0,\quad
t>0,\,x\in \Omega,\nonumber\\
u|_{\partial \Omega} & =0,\quad t>0,\,x\in \partial\Omega, \label{problembdddom}\\
u|_{t=0} & =u_0,\quad x\in \partial\Omega.\nonumber
\end{align}
Let $\{\phi_n\}_{n=1}^\infty
\subset \oH^1_2(\Omega):=\overline{C_0^\infty(\Omega)}\,{}^{H^1_2(\Omega)}$ be an orthonormal basis of $L_2(\Omega)$ consisting of eigenfunctions of the negative
Dirichlet Laplacian with eigenvalues $\lambda_n>0$,
$n\in \iN$, and denote by $\lambda_1$ the smallest such eigenvalue.
Then the solution $u$ of (\ref{problembdddom}) can
be represented via Fourier series as
\begin{equation} \label{uformel}
u(t,x)=\sum_{n=1}^\infty s_{\lambda_n}(t)\,(u_0|\phi_n)\phi_n(x),\quad t\ge 0,\,x\in \Omega,
\end{equation}
where $(\cdot|\cdot)$ stands for the standard inner product in $L_2(\Omega)$, cf.\ \cite[Section 1]{VZ}, the special case $k=g_{1-\alpha}$ can be also found in \cite[Theorem 4.1]{NSY}. By Parseval's identity and since $\partial_\mu s_\mu\le 0$, it follows from (\ref{uformel}) that
\begin{align*}
|u(t,\cdot)|_{L_2(\Omega)}^2 & =\sum_{n=1}^\infty s_{\lambda_n}^2(t)\,|(u_0|\phi_n)|^2\\
& \le s_{\lambda_1}^2(t) \sum_{n=1}^\infty|(u_0|\phi_n)|^2\\
& = s_{\lambda_1}^2(t) |u_0|_{L_2(\Omega)}^2,
\end{align*}
and thus
\begin{equation} \label{udecay}
|u(t,\cdot)|_{L_2(\Omega)}\le s_{\lambda_1}(t)|u_0|_{L_2(\Omega)},\quad t\ge 0,
\end{equation}
cf.\ \cite{VZ}. This decay estimate is optimal as the example $u_0=\phi_1$ with solution $u(t,x)
=s_{\lambda_1}(t)\phi_1(x)$ shows. By means of energy methods (\ref{udecay}) can be generalized to
problems with a uniformly elliptic operator in divergence form with rough coefficients, cf.\ \cite[Corollary 1.1]{VZ}.
(\ref{udecay}) can be further extended to
$r\in (1,\infty)$,
in fact assuming $u_0\in L_r(\Omega)$ and setting $\rho(r):=4(r-1)/r^2$ we have
\[
|u(t,\cdot)|_r\le s_{\lambda_1 \rho(r)}(t)|u_0|_r,\quad t>0,
\]
see \cite[Remark 5.1]{VZ}. We see that the relaxation function $s_\mu(t)$ (with some fixed $\mu>0$) determines the
rate of decay as $t\to \infty$. For example, if $k=g_{1-\alpha}$ with $\alpha\in (0,1)$ it follows from (\ref{smubounds}) that $s_\mu(t)$ decays like $c t^{-\alpha}$. This justifies among others earlier remarks on (\ref{L2intro}) concerning
the case $d>4$.

\smallbreak

{\bf Energy methods for weak solutions.} A further goal of this paper is to prove decay estimates for the time-fractional diffusion equation with the Laplacian being replaced by a more general elliptic operator in divergence form with rough coefficients. More precisely, we want to study
the problem
\begin{align}
\partial_t^\alpha (u-u_0)-\mbox{div}\,\big(A(t,x) \nabla u\big) & =0,\quad t>0,\,x\in \iR^d, \label{IAfracdiff}\\
u|_{t=0} & = u_0,\quad x\in \iR^d. \label{IAfracdiffini}
\end{align}
Here, we assume $\alpha\in (0,1)$, $u_0\in L_1(\iR^d) \cap L_2(\iR^d)$, $A\in L_{\infty, loc}([0,\infty)\times \iR^d;
\iR^{d\times d})$, and that $A$ satisfies a uniform parabolicity condition. We prove that if $u$ is a suitably defined weak solution of (\ref{IAfracdiff}), (\ref{IAfracdiffini}) satisfying appropriate integrability conditions at $x=\infty$,
there holds for all $d\in \iN$
\begin{equation} \label{IthmEDecay1}
|u(t)|_{2} \lesssim t^{-\frac{\alpha d}{d+4}},\quad t>0,
\end{equation}
see Theorem \ref{energydecay}. The basic idea of the proof is to show that for some constant $\mu>0$
\[
\partial_t^\alpha \big(|u|_2-|u_0|_2\big)(t)+\mu |u(t)|_2^{1+\frac{4}{d}}\le 0,\quad t>0,
\]
in the weak sense. This can be achieved by means of Nash's inequality and the so-called {\em $L_p$-norm inequality} for operators of the form $\partial_t(k\ast\cdot)$, which has been established recently in \cite{VZ}, cf.\ Lemma \ref{LemmaL2IN} below. Note that the decay rate in (\ref{IthmEDecay1}) is less than the one we find in the case of the Laplacian (cf.\ (\ref{L2intro})). This phenomenon of a smaller decay rate in the weak setting with rough coefficients does
not occur in the case $\alpha=1$, where the same strategy of proof leads to the optimal decay rate $t^ {-d/4}$. It is an
interesting open problem whether (\ref{IthmEDecay1}) provides the optimal decay rate in the variational setting.

\medbreak

\noindent The paper is organized as follows. In Section 2 we provide some background on the fundamental solution $Z$, including their construction via subordination. Section 3 is devoted to the time-fractional diffusion case. We study the decay properties of $Z$ and $Z(t,\cdot)\star u_0$, respectively. Sections 4 and 5 deal with the general case. We first prove $L_2$-decay estimates for $Z(t,\cdot)\star u_0$ using Plancherel's theorem
and discuss some specific examples of pairs $(k,l)\in {\cal PC}$ (Section 4). Then we establish $L_r$-estimates and
look at the critical dimension case (Section 5). Finally, in Section 6 we discuss decay estimates in the variational setting.
\section{The fundamental solution $Z$} \label{SecFS}
Suppose $(k,l)\in {\cal PC}$. We will describe how the fundamental solution $Z$ to the subdiffusion problem (\ref{subdiffintro}), (\ref{subdiffinitintro}) can be constructed from the Gaussian heat kernel $H$ by means of the subordination principle for abstract Volterra equations with completely positive kernels, cf.\ Pr\"uss \cite[Chapter 4]{JanI}. See also the proof of \cite[Theorem 3]{Koch11}. We will also show that $Z(t,\cdot)$ is a probability density function on $\iR^d$ for all $t>0$.

By definition, the fundamental solution $Z(t,x)$ to the subdiffusion problem (\ref{subdiffintro}), (\ref{subdiffinitintro}) is a distributional solution of
\[
\partial_t\big(k\ast [Z-Z_0]\big)-\Delta Z=0,\;\;t>0,\,x\in \iR^d,\quad Z|_{t=0}=Z_0:=\delta_0,\,x\in \iR^d,
\]
where $\delta_0$ stands for the Dirac delta distribution. Throughout this paper the Fourier transform of $v\in \mathcal{S}(\iR^d)$ is defined by
\begin{equation} \label{FourierDef}
\tilde{v}(\xi)=\int_{\iR^d} e^{-i x\cdot \xi}v(x)\,dx,
\end{equation}
extended as usual to $\mathcal{S}'(\iR^d)$. Taking the Fourier transform w.r.t.\ $x$ in the problem for $Z$ we obtain
(cf.\ Section 1)
\[
\tilde{Z}(t,\xi)=s(t,|\xi|^2),\quad t\ge 0,\,\xi \in \iR^d,
\]
where the relaxation function $s(t,\mu)$ is defined via the Volterra equation (\ref{smudef}). In what follows we will construct a function $Z$ that enjoys the latter property.

Note first that $(k,l)\in {\cal PC}$ implies that the kernel $l$ is completely positive, see \cite{CN} and \cite{JanI}.
Denoting by $\hat{f}$ the Laplace transform of $f:\iR_+\rightarrow \iR$, this in turn implies that $\varphi(\lambda)=\lambda\hat{k}(\lambda)$, $\lambda>0$, is a Bernstein function and that for every $\tau\ge 0$, the function $\psi_\tau(\lambda)=\exp(-\tau \varphi(\lambda))$,
$\lambda>0$, is completely monotone, by \cite[Proposition 4.5]{JanI}. Further, $\psi_\tau(\lambda)$ is bounded by
$e^{-\tau \varphi(0+)}$. Note that $\varphi(\lambda)=1/\hat{l}(\lambda)$, and thus $\varphi(0+)=1/|l|_{L_1(\iR_+)}$, which is $0$ if and only if $l\notin L_1(\iR_+)$.
By Bernstein's theorem (see e.g.\ \cite[Section 4.1]{JanI} or \cite[Theorem 1.4]{SSV10}) there exist unique nondecreasing functions $w(\cdot,\tau)\in BV(\iR_+)$
normalized by $w(0,\tau)=0$ and left-continuity such that
\[
\hat{w}(\lambda,\tau)=\int_0^\infty e^{-\lambda \sigma}w(\sigma,\tau)\,d\sigma=\,\frac{\psi_\tau(\lambda)}
{\lambda},\quad \lambda>0.
\]
The function $w(t,\tau)$ is the so-called {\em propagation function} associated with the completely positive kernel $l$,
cf.\ \cite[Section 4.5]{JanI}. Some important properties of $w(t,\tau)$ can be found in \cite[Proposition 4.9]{JanI}.
Among others, $w(\cdot,\cdot)$ is Borel measurable on $\iR_+\times \iR_+$, $w(t,\cdot)$ is nonincreasing and right-continuous on $\iR_+$, and $w(t,0)=w(t,0+)=1$ as well as
$w(t,\infty)=0$ for all $t>0$. Moreover, the relaxation function $s(t,\mu)$ is represented by
\begin{equation} \label{smurepresent}
s(t,\mu)=-\int_0^\infty e^{-\mu \tau}\,d_\tau w(t,\tau),\quad t>0,\,\mu\ge 0,
\end{equation}
in particular $-\int_0^\infty d_\tau w(t,\tau)=1$.

We now set
\[
Z(t,x)=-\int_0^\infty H(\tau,x)\,d_\tau w(t,\tau),\quad t>0,\,x\in \iR^d.
\]
Then $Z(t,x)$ is nonnegative, and $|Z(t,\cdot)|_1=1$ for all $t>0$, since $H$ enjoys these properties.
Taking the Fourier transform w.r.t.\ $x$ and using (\ref{smurepresent}) we obtain
\[
\tilde{Z}(t,\xi)=-\int_0^\infty \tilde{H}(t,\xi)\,d_\tau w(t,\tau)=-\int_0^\infty e^{-\tau |\xi|^2}\,d_\tau w(t,\tau)
=s(t,|\xi|^2).
\]
Thus $Z$ is the desired fundamental solution to (\ref{subdiffintro}).

The following lemma provides a formula for the mean square displacement. The basic idea of the proof can be found
already in \cite[p.\ 19]{Metz} (time-fractional case with $d=1$) and \cite[p.\ 268]{Koch08} (ultraslow diffusion with arbitrary $d$).
\begin{lemma} \label{MSD}
Let $(k,l)\in {\cal PC}$ and $Z$ be the fundamental solution to the diffusion problem (\ref{subdiffintro}), (\ref{subdiffinitintro}). Let $m(t)$ be defined as in (\ref{MSDdef}). Then
\[
m(t)= \,2d\,(1\ast l)(t),\quad t> 0.
\]
\end{lemma}
{\em Proof.} For the Laplace transform of $m$ we have
\begin{align*}
\hat{m}(\lambda) & =\int_0^\infty \int_{\iR^d} |x|^2 e^{-\lambda t}Z(t,x)\,dx\,dt\\
& =\int_0^\infty \int_{\iR^d} \big(-\Delta_\xi e^{-ix\cdot \xi}\big)|_{\xi=0} e^{-\lambda t}Z(t,x)\,dx\,dt\\
& = \Big( -\Delta_\xi \int_0^\infty \int_{\iR^d} e^{-ix\cdot \xi}e^{-\lambda t}Z(t,x)\,dx\,dt\Big)\Big|_{\xi=0}\\
& = \Big( -\Delta_\xi \hat{\tilde{Z}}(\lambda,\xi) \Big)\Big|_{\xi=0}\\
& = \Big( -\Delta_\xi \hat{s}(\lambda,|\xi|^2) \Big)\Big|_{\xi=0}\\
& = \Big( -\Delta_\xi \big[\frac{1}{\lambda(1+|\xi|^2\hat{l}(\lambda))}\big] \Big)\Big|_{\xi=0}\\
& = 2d\,\frac{\hat{l}(\lambda)}{\lambda},
\end{align*}
and thus the claim follows by inversion of the Laplace transform. \hfill $\square$

${}$

Let us illustrate Lemma \ref{MSD} by looking at the time-fractional diffusion case, where $k=g_{1-\alpha}$ and
$l=g_\alpha$ with some $\alpha\in (0,1)$. The formula for $m(t)$ gives
\[
m(t)=2d (1\ast g_\alpha)(t)=2d g_{1+\alpha}(t)=\,\frac{2d}{\Gamma(1+\alpha)}\,t^\alpha,\quad t>0,
\]
which is in accordance with our remarks in the paragraph following (\ref{MSDdef}) in Section 1.
\section{Time-fractional diffusion} \label{fractionalheat}
In this section we study the subdiffusion problem (\ref{subdiffintro}),
(\ref{subdiffinitintro}) in the important special case $k=g_{1-\alpha}$ with $\alpha\in (0,1)$. That is, we consider the problem
\begin{align}
\partial_t^\alpha (u-u_0)-\Delta u & =0,\quad t>0,\,x\in \iR^d, \label{fracdiff}\\
u|_{t=0} & = u_0,\quad x\in \iR^d. \label{fracdiffini}
\end{align}
Under appropriate conditions on $u_0$ the solution of (\ref{fracdiff}), (\ref{fracdiffini}) can be represented as
\begin{equation} \label{solformulaTF}
u(t,x)=\int_{\iR^d} Z(t,x-y)u_0(y)\,dy,
\end{equation}
where $Z$ is the fundamental solution corresponding to (\ref{fracdiff}),(\ref{fracdiffini}),
see \cite{Koch}. It is known (see e.g.\ \cite{Koch90}, \cite{SchnWyss}) that
\[
Z(t,x)=\pi^{-\frac{d}{2}}t^{\alpha-1}|x|^{-d} H^{20}_{12}\big(\frac{1}{4}|x|^2t^{-\alpha}
\big|^{(\alpha,\alpha)}_{(d/2,1), (1,1)}\big),\quad t>0,\,x\in \iR^d\setminus\{0\},
\]
where $H$ denotes the Fox $H$-function (\cite{KiSa}, \cite{KST}). As the $H$-function is a rather complicated object, this
representation of $Z$ is not so useful for deriving estimates for $Z$ directly. However, using the analytic and asymptotic properties of $H$, one can obtain the subsequent sharp estimates, which can be found in \cite{Koch}, see also \cite{Koch90}.
\begin{prop} \label{Kochubei1}
Set $R=t^{-\alpha}|x|^2$. Then
\begin{align}
Z(t,x) & \le C t^{-\frac{\alpha d}{2}} \exp\left(-\sigma R^{\frac{1}{2-\alpha}}\right), \quad \mbox{if}\;R\ge 1, \label{Z1}\\
Z(t,x) & \le C t^{-\alpha} |x|^{-d+2},  \quad\quad\quad \mbox{if}\;R\le 1\;\mbox{and}\;d\ge 3,\label{Z2}\\
Z(t,x) & \le C t^{-\alpha} \left(|\log R|+1\right),  \quad \mbox{if}\;R\le 1\;\mbox{and}\;d=2,\label{Z3}\\
Z(t,x) & \le C t^{-\frac{\alpha}{2}}, \quad\quad\quad \mbox{if}\;R\le 1\;\mbox{and}\;d=1,\label{Z4}\\
|\nabla Z(t,x)| & \le C t^{-\frac{\alpha (d+1)}{2}} \exp\left(-\sigma R^{\frac{1}{2-\alpha}}\right)   , \quad \mbox{if}\;R\ge 1, \nonumber\\
|\nabla Z(t,x)| & \le C  t^{-\alpha} |x|^{-d+1}  , \quad \mbox{if}\;R\le 1\;\mbox{and}\;d\ge 2, \nonumber\\
|\nabla Z(t,x)| & \le C  t^{-\alpha} , \quad \quad\quad\mbox{if}\;R\le 1\;\mbox{and}\;d=1. \nonumber
\end{align}
Here $C=C(\alpha,d)$ is a positive constant which may differ from line to line.
\end{prop}
We point out that for $d\ge 2$ the fundamental solution $Z(t,x)$ has a singularity at the origin $x=0$, which is
a fundamental difference to the heat kernel and which leads to a restriction concerning $p$-integrability on $\iR^d$.
\subsection{$L_p(\iR^d)$-estimates for $Z$ and the solution}
Proposition \ref{Kochubei1} allows us to estimate the $L_p(\iR^d)$-norm of $Z(t,\cdot)$ for $t>0$.
We decompose the corresponding integral as follows.
\[
|Z(t)|_p^p \le \int_{\{R\ge 1\}}Z(t,x)^p\,dx+\int_{\{R\le 1\}}Z(t,x)^p\,dx.
\]
For all dimensions $d$ and $1<p<\infty$, we have in view of (\ref{Z1})
\begin{align*}
\int_{\{R\ge 1\}}Z(t,x)^p\,dx  & \lesssim \int_{\{R\ge 1\}} t^{-\frac{\alpha d p}{2}} \exp\left(-\sigma p R^{\frac{1}{2-\alpha}}\right)\,dx \\
& \lesssim \int_{t^\frac{\alpha}{2}}^\infty  t^{-\frac{\alpha d p}{2}} \exp\left(-\sigma p
 \big(r^2 t^{-\alpha}\big)^{\frac{1}{2-\alpha}}\right)\,r^{d-1}\,dr\\
& \lesssim t^{ -\frac{\alpha d p}{2}+\frac{\alpha d}{2}} \int_1^\infty \exp\left(-\sigma p\, s^\frac{2}{2-\alpha}\right)\,s^{d-1}\,ds,
\end{align*}
and thus
\begin{equation} \label{goodpart}
\left(\int_{\{R\ge 1\}}Z(t,x)^p\,dx\right)^\frac{1}{p} \lesssim t^{ -\frac{\alpha d }{2}\,\left(1-\frac{1}{p}\right)},\quad t>0,\quad \mbox{for all}\;1<p<\infty.
\end{equation}

We come now to the estimate for the integral where $R\le 1$. In case $d=1$ we have in view of (\ref{Z4}) for all $1<p<\infty$
\begin{align*}
\int_{\{R\le 1\}}Z(t,x)^p\,dx & \lesssim \int_{\{R\le 1\}} t^{-\frac{\alpha p}{2}}\,dx  \lesssim
\int_0^{t^\frac{\alpha}{2}} t^{-\frac{\alpha p}{2}}\,dx \lesssim t^{-\frac{\alpha p}{2}+\frac{\alpha}{2}}.
\end{align*}
If $d=2$ we may estimate as follows, employing (\ref{Z3}).
\begin{align*}
\int_{\{R\le 1\}}Z(t,x)^p\,dx & \lesssim \int_{\{R\le 1\}} t^{-\alpha p} \left(|\log R|+1\right)^p\,dx\\
& \lesssim \int_0^{t^\frac{\alpha}{2}} t^{-\alpha p} \left(|\log (r^2 t^{-\alpha})|+1\right)^p\,r\,dr\\
& \lesssim \int_0^1 t^{-\alpha p+\alpha} \left(|\log (s^2)|+1\right)^p\,s\,ds \lesssim   t^{-\alpha p+\alpha},
\end{align*}
for all $1<p<\infty$. For $d\ge 3$ we use (\ref{Z2}) to obtain
\begin{align*}
\int_{\{R\le 1\}}Z(t,x)^p\,dx & \lesssim \int_{\{R\le 1\}} t^{-\alpha p} |x|^{(-d+2)p}\,dx
\lesssim \int_0^{t^\frac{\alpha}{2}} t^{-\alpha p} r^{(-d+2)p} r^{d-1}\,dr \\
& \lesssim   t^{-\alpha p+\frac{\alpha}{2}(d+[2-d]p)} \int_0^1 s^{(-d+2)p} s^{d-1}\,ds
\lesssim t^{-\frac{\alpha d}{2}\, (p-1)},
\end{align*}
whenever the last integral is finite, that is, whenever
\[
p<\,\frac{d}{d-2}=:\kappa(d).
\]
Setting $\kappa(1)=\kappa(2)=\infty$ and
combining the previous estimates we see that
\[
\left(\int_{\{R\le 1\}}Z(t,x)^p\,dx\right)^\frac{1}{p} \lesssim t^{ -\frac{\alpha d }{2}\,\left(1-\frac{1}{p}\right)},\quad t>0,\quad \mbox{for all}\;1<p<\kappa(d).
\]

Observe also that $Z(t)\in L_\infty(\iR)$, $t>0$, provided that $d=1$, and we have the estimate
$|Z(t)|_\infty \lesssim t^{-\alpha/2}$.

Summarizing we have proved
\begin{satz} \label{ZLpnorm}
Let $d\in\iN$ and $\kappa(d)$ be as above. Then $Z(t)$ belongs to $L_p(\iR^d)$ for all $t>0$, provided that
$1\le p<\kappa(d)$, and there holds
\begin{equation} \label{ZLpEst}
|Z(t)|_p  \lesssim t^{ -\frac{\alpha d }{2}\,\left(1-\frac{1}{p}\right)},\quad t>0.
\end{equation}
(\ref{ZLpEst}) remains true for $d=1$ and $p=\kappa(1)=\infty$.
\end{satz}
\begin{bemerk1}
{\em
 $Z(t)$ fails to belong to $L_p(\iR^d)$ for $d\ge 4$ and
$p\ge \kappa(d)=d/(d-2)$.
In fact, taking the Fourier transform w.r.t.\ the spatial variables we see that the Fourier transform of $Z(t,x)$ solves the fractional differential equation
\[
\partial_t^\alpha (\tilde{Z}-1)+|\xi|^2 \tilde{Z}=0,\;\; t>0,\quad\quad \tilde{Z}(0)=1.
\]
Thus
\[
\tilde{Z}(t,\xi)=s(t,|\xi|^2) =E_\alpha\left(-|\xi|^2 t^\alpha\right),
\]
where $E_\alpha(z)$ denotes the Mittag-Leffler function, defined by
\[
E_\alpha(z)=\sum_{j=0}^\infty \,\frac{z^j}{\Gamma(\alpha j+1)}\,,\quad z\in\iC.
\]
Employing the bounds from (\ref{smubounds}), a short computation shows that the Mittag-Leffler function satisfies the estimate
\begin{equation} \label{ML}
\frac{1}{1 + \Gamma (1-\alpha) x} \leq E_{\alpha}(-x) \leq \frac{1}{1+\frac{x}{\Gamma(1+\alpha)}},\,\quad x\ge 0,
\end{equation}
see also \cite[Example 6.1]{VZ}. An upper bound of the form $E_\alpha(-x)\le C(\alpha)/(1+x)$, $x\ge 0$, can also be found in \cite{MNV}.

Suppose now that $Z(t)$ belongs to $L_p(\iR^d)$ with some $p\le 2$. Then the Hausdorff-Young inequality
implies that
\[
|\tilde{Z}(t)|_{p'}\lesssim |Z(t)|_p <\infty,
\]
where $p'$ is the conjugate exponent of $p$. From this and the representation of $\tilde{Z}$ as well as
(\ref{ML}) it follows that
\[
\int_{\iR^d} \,\frac{d\xi}{(1+|\xi|^2 t^\alpha)^{p'}}\,<\infty,
\]
which in turn implies (by changing to polar coordinates) that
\[
\int_0^\infty \,\frac{r^{d-1}}{(1+r^2)^{p'}}\,dr<\infty.
\]
Hence $2p'-(d-1)>1$, which is equivalent to $p<d/(d-2)$. \hfill $\square$
}
\end{bemerk1}

${}$

\noindent We next examine the critical case $p=\frac{d}{d-2}$ for $d\ge 3$.
We know that $Z(t)$ does not belong to $L_p(\iR^d)$. However, it lies in the corresponding weak $L_p$-space as the following theorem shows. This observation will be crucial, among others, to obtain optimal
decay estimates for $|u(t)|_2$ for $d\ge 4$.
\begin{satz} \label{ZweakLp}
Let $d\ge 3$ and $t>0$. Then $Z(t)$ belongs to the space $L_{\frac{d}{d-2},\,\infty}(\iR^d)$, and there holds the estimate
\[
|Z(t)|_{\frac{d}{d-2},\,\infty} \lesssim \,t^{-\alpha},\quad t>0.
\]
\end{satz}
{\em Proof.} Set $p=\frac{d}{d-2}$. We need to estimate
\[
|Z(t)|_{p,\,\infty}=\sup \left\{\lambda\, d_{Z(t)}(\lambda)^{\frac{1}{p}}:\,\lambda>0\right\},
\]
where
\[
d_{Z(t)}(\lambda)=|\{x\in \iR^d:\,Z(t,x)>\lambda\}|
\]
denotes the distribution function of $Z(t)$. Using again the similarity variable $R=t^{-\alpha}|x|^2$ we have
\begin{align*}
|Z(t)|_{p,\,\infty} & \le 2\left(|Z(t)\chi_{\{R\le 1\}}(t)|_{p,\,\infty}+
|Z(t)\chi_{\{R\ge 1\}}(t)|_{p,\,\infty}\right).
\end{align*}
Employing (\ref{goodpart}), we find that
\[
|Z(t)\chi_{\{R\ge 1\}}(t)|_{p,\,\infty}
\le |Z(t)\chi_{\{R\ge 1\}}(t)|_{p}
\le C t^{ -\frac{\alpha d }{2}\,\left(1-\frac{1}{p}\right)}
= C t^{-\alpha}.
\]
For the term with $R\le 1$ we use (\ref{Z2}) to estimate as follows.
\begin{align*}
d_{Z(t)\chi_{\{R\le 1\}}(t)}(\lambda) & =|\{x\in \iR^d:\,Z(t,x)>\lambda\;\mbox{and}\;R\le 1\}|\\
& \le |\{x\in \iR^d:\,\lambda< C t^{-\alpha} |x|^{-d+2}\}|\\
& = |\{x\in \iR^d:\,|x|< \left(C t^{-\alpha}\lambda^{-1}\right)^{\frac{1}{d-2}}\}|\\
& \le C_1 \left( t^{-\alpha}\lambda^{-1}\right)^{\frac{d}{d-2}}.
\end{align*}
This shows that
\[ d_{Z(t)\chi_{\{R\le 1\}}(t)}(\lambda)^{1/p}\le C_1^{1/p} t^{-\alpha}\lambda^{-1},
\]
and thus
\[
|Z(t)\chi_{\{R\le 1\}}(t)|_{p,\,\infty}\lesssim\,t^{-\alpha}.
\]
This proves the theorem. \hfill $\square$

${}$

\noindent We come now to decay estimates for the $L_r(\iR^d)$-norm of the solution $u$ of
(\ref{fracdiff}), (\ref{fracdiffini}) given by formula (\ref{solformulaTF}).
\begin{satz} \label{Lrestimate}
(i) Let $d\in \iN$, $1\le p<\kappa(d)$, $1\le q,r \le \infty$, $1+\frac{1}{r}=\frac{1}{p}+
\frac{1}{q}$, and $u_0\in L_q(\iR^d)$. Then for $u(t)=Z(t)\star u_0$ we have
\[
|u(t)|_r \lesssim   t^{ -\frac{\alpha d }{2}\,\left(1-\frac{1}{p}\right)},\quad t>0.
\]
This estimate remains true for $d=q=1$ and $p=r=\infty$.

\smallbreak

\noindent (ii) Let $d\ge 3$, $1< q,r < \infty$, $\frac{1}{r}+\frac{2}{d}=\frac{1}{q}$, and $u_0\in L_q(\iR^d)$. Then for $u(t)=Z(t)\star u_0$ we have
\[
|u(t)|_r \lesssim   t^{ -\alpha},\quad t>0.
\]
(iii) Let $d\ge 3$ and $u_0\in L_1(\iR^d)$. Then for $u(t)=Z(t)\star u_0$ we have
\[
|u(t)|_{\frac{d}{d-2},\infty}\lesssim   t^{ -\alpha},\quad t>0.
\]
\end{satz}
{\em Proof.} (i) By Young's inequality and Theorem \ref{ZLpnorm} we have
\[
|u(t)|_r \le |Z(t)|_p |u_0|_q \lesssim  t^{ -\frac{\alpha d }{2}\,\left(1-\frac{1}{p}\right)}, \quad t>0.
\]

(ii) We apply Young's inequality for weak type spaces (cf.\ \cite[Theorem 1.4.24]{Graf}) with $p=\kappa(d)=\frac{d}{d-2}$
and invoke Theorem \ref{ZweakLp}. This yields
\[
|u(t)|_r \le C(p,q,r) |Z(t)|_{p,\infty} |u_0|_q \lesssim   t^{ -\alpha},\quad t>0.
\]

(iii) The assertion follows from Theorem \ref{ZweakLp} and Young's inequality for weak type spaces (see \cite[Theorem 1.2.13]{Graf}).
\hfill $\square$
\begin{korollar} \label{Lrkorollar}
Let $1<r<\infty$, $d\in \iN$, $u_0\in L_1(\iR^d)\cap L_r(\iR^d)$, and $u(t)=Z(t)\star u_0$. Then
\begin{align*}
|u(t)|_r \,& \lesssim\, t^{ -\frac{\alpha d }{2}\,\left(1-\frac{1}{r}\right)},\quad t>0,\;\; \mbox{if}\;\;d<3,\\
|u(t)|_r \,& \lesssim\, t^{ -\frac{\alpha d }{2}\,\left(1-\frac{1}{r}\right)},\quad t>0,\;\; \mbox{if}\;\;d\ge 3\;\mbox{and}\;d<\frac{2r}{r-1},\\
|u(t)|_{r,\infty} \,& \lesssim\,t^{-\alpha},\quad t>0,\;\; \mbox{if}\;\;d\ge 3\;\mbox{and}\;d=\frac{2r}{r-1},\\
|u(t)|_r \,& \lesssim\,t^{-\alpha},\quad t>0,\;\; \mbox{if}\;\;d\ge 3\;\mbox{and}\;d>\frac{2r}{r-1}.
\end{align*}
\end{korollar}
{\em Proof.} If $d<3$ or we have both $d\ge 3$ and $d<\frac{2r}{r-1}$ then $r<\kappa(d)$ and the first two
estimates follow from Theorem \ref{Lrestimate} (i) with $q=1$. The third estimate is a consequence of
Theorem \ref{Lrestimate} (iii), since $d=\frac{2r}{r-1}$ is equivalent
to $r=\kappa(d)$ whenever $d\ge 3$. To show the last estimate, observe first that the assumptions on $r$ and $d$
imply that there is $q\in (1,r)$ such that $\frac{1}{r}+\frac{2}{d}=\frac{1}{q}$. The assertion then follows from
Theorem \ref{Lrestimate} (ii). \hfill $\square$

${}$

\noindent Specializing Corollary \ref{Lrkorollar} to $r=2$ we obtain
\begin{korollar} \label{L2estimate}
Let $d\in \iN$ and $u_0\in L_1(\iR^d)\cap L_2(\iR^d)$ and $u(t)=Z(t)\star u_0$. Then
\begin{align}
|u(t)|_2\,& \lesssim\,t^{-\min\{\frac{\alpha d}{4},\alpha\}},\quad t>0,\;d\in \iN\setminus\{4\}, \label{L2decay}\\
|u(t)|_{2,\infty} & \lesssim\,t^{-\alpha},\quad t>0,\;d=4.\nonumber
\end{align}
\end{korollar}
Corollary \ref{Lrkorollar} exhibits the critical dimension phenomenon discussed already in Section 1. The critical dimension (which is in general not an integer) is here given by $d_{crit}=2r/(r-1)$. As long as the dimension $d$ is below $d_{crit}$
(and at least $3$) the decay rate increases with the dimension, whereas for any $d>d_{crit}$ the decay rate is the same,
namely $t^{ -\alpha}$, which coincides with the decay rate for the corresponding problem on a bounded domain with homogeneous Dirichlet boundary condition.
\subsection{Gradient estimates}
We turn now to $L_p$-estimates for the spatial gradient of $Z(t,x)$. For $d\in \iN$ we set $\kappa_1(d):=\frac{d}{d-1}$ for $d\ge 2$ and $\kappa(1):=\infty$
\begin{satz} \label{GradZ}
(i) Let $d\in \iN$. Then $\nabla Z(t)$ belongs to $L_p(\iR^d;\iR^d)$ for all $t>0$, provided that
$1\le p<\kappa_1(d)$, and there holds
\begin{equation} \label{gradZdecay}
|\nabla Z(t)|_p \lesssim t^{ -\frac{\alpha}{2} -\frac{\alpha d }{2}\,\left(1-\frac{1}{p}\right)},\quad t>0.
\end{equation}
(\ref{gradZdecay}) remains valid for $d=1$ and $p=\kappa_1(d)=\infty$.

(ii) Let $d\ge 2$ and $t>0$. Then $\nabla Z(t)$ belongs to $L_{\frac{d}{d-1},\,\infty}(\iR^d;\iR^d)$,
and we have
\[
|\nabla Z(t)|_{\frac{d}{d-1},\,\infty} \lesssim \,t^{-\alpha},\quad t>0.
\]
\end{satz}
{\em Proof.} The proof uses the gradient estimates for $Z$ from Proposition \ref{Kochubei1} and is analogous to that of Theorem \ref{ZLpnorm} and Theorem \ref{ZweakLp}, respectively. \hfill $\square$

${}$

\noindent By means of the $L_p(\iR^d)$-estimates for $\nabla Z(t,x)$ from Theorem \ref{GradZ} we can derive temporal decay estimates for the gradient of the solution. We record these estimates in the following theorem.
Using $\nabla u(t,\cdot)=\nabla Z(t,\cdot)\star u_0$ the proof is analogous
to the one of Theorem \ref{Lrestimate}.
\begin{satz} \label{gradLrestimate}
(i) Let $d\in \iN$, $1\le p<\kappa_1(d)$, $1\le q,r \le \infty$, $1+\frac{1}{r}=\frac{1}{p}+
\frac{1}{q}$, and $u_0\in L_q(\iR^d)$. Then for $u(t)=Z(t)\star u_0$ we have
\[
|\nabla u(t)|_r \lesssim   t^{ -\frac{\alpha}{2} -\frac{\alpha d }{2}\,\left(1-\frac{1}{p}\right)},\quad t>0.
\]
(ii) Let $d\ge 2$, $1< q,r < \infty$, $\frac{1}{r}+\frac{2}{d}=\frac{1}{q}$, and $u_0\in L_q(\iR^d)$. Then for $u(t)=Z(t)\star u_0$ we have
\[
|\nabla u(t)|_r \lesssim   t^{ -\alpha},\quad t>0.
\]
(iii) Let $d\ge 2$ and $u_0\in L_1(\iR^d)$. Then for $u(t)=Z(t)\star u_0$ we have
\[
|\nabla u(t)|_{\frac{d}{d-1},\infty}\lesssim   t^{ -\alpha},\quad t>0.
\]
\end{satz}

${}$

\noindent Arguing analogously to the proof of Corollary \ref{Lrkorollar} we obtain
\begin{korollar} \label{GradLrkorollar}
Let $1<r<\infty$, $d\in \iN$, $u_0\in L_1(\iR^d)\cap L_r(\iR^d)$, and $u(t)=Z(t)\star u_0$. Then
\begin{align*}
|\nabla u(t)|_r \,& \lesssim\, t^{ -\frac{\alpha}{2} -\frac{\alpha d }{2}\,\left(1-\frac{1}{r}\right)},\quad t>0,\;\; \mbox{if}\;\;d=1,\\
|\nabla u(t)|_r \,& \lesssim\, t^{ -\frac{\alpha}{2} -\frac{\alpha d }{2}\,\left(1-\frac{1}{r}\right)},\quad t>0,\;\; \mbox{if}\;\;d\ge 2\;\mbox{and}\;d<\frac{
r}{r-1},\\
|\nabla u(t)|_{r,\infty} \,& \lesssim\,t^{-\alpha},\quad t>0,\;\; \mbox{if}\;\;d\ge 2\;\mbox{and}\;d=\frac{r}{r-1},\\
|\nabla u(t)|_r \,& \lesssim\,t^{-\alpha},\quad t>0,\;\; \mbox{if}\;\;d\ge 2\;\mbox{and}\;d>\frac{r}{r-1}.
\end{align*}
\end{korollar}
In the important special case $r=2$ the picture is as follows.
\begin{korollar} \label{GradL2estimate}
Let $d\in \iN$ and $u_0\in L_1(\iR^d)\cap L_2(\iR^d)$ and $u(t)=Z(t)\star u_0$. Then
\begin{align}
|\nabla u(t)|_2\,& \lesssim\,t^{-\min\{\frac{\alpha}{2}+\frac{\alpha d}{4},\alpha\}},\quad t>0,\;d\in \iN\setminus\{2\}, \nonumber\\
|\nabla u(t)|_{2,\infty} & \lesssim\,t^{-\alpha},\quad t>0,\;d=2.\nonumber
\end{align}
\end{korollar}
Corollary \ref{GradLrkorollar} shows a critical dimension phenomenon for the gradient estimates with critical dimension $d_{grad,crit}=r/(r-1)<d_{crit}$.
\subsection{Large time behaviour of $Z(t,\cdot)\star u_0$}
In this subsection, we want to show that for integrable initial data $u_0$ the asymptotic behaviour of
$Z(t,\cdot)\star u_0$ as $t\to\infty$ is described by a multiple of $Z(t,x)$. The corresponding result for the
heat equation is well-known, see e.g.\ \cite[Prop.\ 48.6]{QS}.

For our strategy of proof we need the following decomposition lemma from \cite{DZua}.
\begin{lemma} \label{decomp}
Suppose $f\in L_1(\iR^d)$ such that $\int_{\iR^d} |x|\,|f(x)|\,dx<\infty$. Then there exists $F\in L_1(\iR^d;\iR^d)$ such that
\[
f=\left(\int_{\iR^d} f(x)\,dx\right)\,\delta_0+\mbox{{\em div}}\,F
\]
in the distributional sense and
\[
|F|_{L_1(\iR^d;\iR^d)}\le C_d \int_{\iR^d} |x|\,|f(x)|\,dx.
\]
\end{lemma}
We have now the following result.
\begin{satz} Let $d\in \iN$ and $1\le p<\kappa_1(d)$. Let $u_0\in L_1(\iR^d)$ and set $M=\int_{\iR^d} u_0(y)\,dy$.

\noindent (i) There holds
\[
t^{\frac{\alpha d }{2}\,\left(1-\frac{1}{p}\right)}|u(t)-MZ(t)|_{p}\rightarrow 0,\quad\mbox{as}\;t\rightarrow \infty.
\]
(ii) Assume in addition that $||x|u_0|_1<\infty$. Then
\[
t^{\frac{\alpha d }{2}\,\left(1-\frac{1}{p}\right)}|u(t)-MZ(t)|_{p}\lesssim t^{ -\frac{\alpha}{2}},\quad t>0.
\]
Moreover, in the limit case $p=\kappa_1(d)$ we have
\[
t^{\frac{\alpha }{2}}|u(t)-MZ(t)|_{\kappa_1(d),\,\infty}\lesssim t^{ -\frac{\alpha}{2}},\quad t>0.
\]
\end{satz}
{\em Proof.} The strategy of the proof is the same as in \cite[p.\ 14, 15]{Zua}.

 (a) Suppose first that $u_0\in L_1(\iR^d)$ is such that $\int_{\iR^d} |x|\,|u_0(x)|\,dx<\infty$.
By Lemma \ref{decomp} there exists $\phi\in L_1(\iR^d;\iR^d)$ such that
\[
u_0=M\delta_0+\mbox{ div}\,\phi
\]
and $|\phi|_1\le C_d||x|u_0|_1$. Consequently,
\begin{align*}
u(t,x)&=M\left(Z(t,\cdot)\star \delta_0\right)(x)+\left(Z(t,\cdot)\star \mbox{ div}\,\phi(\cdot)\right)(x)\\
& =M Z(t,x)+(\nabla Z(t,\cdot)\star \phi)(x),
\end{align*}
which yields
\begin{equation} \label{asym1}
u(t,x)-M Z(t,x)=(\nabla Z(t,\cdot)\star \phi)(x).
\end{equation}
By Young's inequality it follows that for any $1\le p<\kappa_1(d)$
\[
|u(t)-MZ(t)|_{p}\le |\nabla Z(t)|_{p} |\phi|_1 \lesssim  |\nabla Z(t)|_{p} \big||x|u_0\big|_1\lesssim
t^{ -\frac{\alpha}{2} -\frac{\alpha d }{2}\,\left(1-\frac{1}{p}\right)},
\]
where we used Theorem \ref{GradZ}. Hence
\[
t^{\frac{\alpha d }{2}\,\left(1-\frac{1}{p}\right)}|u(t)-MZ(t)|_{p}\lesssim t^{ -\frac{\alpha}{2}},
\]
which is the first part of assertion (ii). The second part follows from
(\ref{asym1}) by applying Young's inequality for weak $L_p$-spaces
(\cite[Theorem 1.2.13]{Graf}).

(b) To prove (i) we choose a sequence $(\eta_j)\subset C_0^\infty(\iR^d)$ such that
$\int_{\iR^d} \eta_j\,dx=M$ for all $j$ and $\eta_j \rightarrow u_0$ in $L_1(\iR^d)$. For each $j$ we have by Part (a) and by Theorem \ref{ZLpnorm}
\begin{align*}
|u(t)-MZ(t)|_{p} & \le |Z(t)\star(u_0-\eta_j)|_p+
  |Z(t)\star \eta_j-MZ(t)|_p\\
  & \le  |Z(t)|_p |u_0-\eta_j|_1+C(j)\,t^{ -\frac{\alpha}{2} -\frac{\alpha d }{2}\,\left(1-\frac{1}{p}\right)}\\
  & \le C_1 t^{ -\frac{\alpha d }{2}\,\left(1-\frac{1}{p}\right)}|u_0-\eta_j|_1+
C(j)\,t^{ -\frac{\alpha}{2} -\frac{\alpha d }{2}\,\left(1-\frac{1}{p}\right)},
\end{align*}
and therefore
\[
t^{ \frac{\alpha d }{2}\,\left(1-\frac{1}{p}\right)}|u(t)-MZ(t)|_{p}\le C_1|u_0-\eta_j|_1+C(j)\,t^{ -\frac{\alpha}{2}},
\]
which implies
\[
\limsup_{t\to \infty}\,t^{ \frac{\alpha d }{2}\,\left(1-\frac{1}{p}\right)}|u(t)-MZ(t)|_{p}\le C_1|u_0-\eta_j|_1.
\]
Assertion (i) follows by sending $j\to \infty$. \hfill $\square$
\section{General subdiffusion equations, $L_2$-estimates} \label{SecL2}
Let $(k,l)\in {\cal PC}$. Assuming that $u_0\in L_1(\iR^d)\cap L_2(\iR^d)$ we want to derive $L_2$-decay estimates for solutions of
\begin{equation} \label{subdiffgen}
\partial_t \big(k\ast [u-u_0]\big)-\Delta u=0,\quad
t>0,\,x\in \iR^d,
\end{equation}
together with
\begin{equation} \label{initcond}
u|_{t=0}=u_0,\,x\in \iR^d,
\end{equation}
 by means of harmonic analysis, the main tool being Plancherel's theorem.

Under appropriate assumptions on the initial value $u_0$ (see
Kochubei \cite{Koch11}) the solution of (\ref{subdiffgen}), (\ref{initcond}) is given by
\begin{equation} \label{gensol}
u(t,x)=\int_{\iR^N} Z(t,x-y)u_0(y)\,dy.
\end{equation}
In the following we will assume that $u$ is defined by (\ref{gensol}).

\begin{satz} \label{theoremlower}
Let $d\in \iN$ and $u_0\in L_1(\iR^d)\cap L_2(\iR^d)$. Assume that $k$ satisfies condition ($\mathcal{PC}$) and that $u$ is given by (\ref{gensol}). Assume further that
$\widetilde{u_0}(0)\neq 0$. Then
\[
|u(t)|_2 \gtrsim k(t)^{\min\{1,\frac{d}{4}\}},\quad \mbox{a.a.}\;t\ge 1.
\]
\end{satz}
{\em Proof.}
Let $\rho_0>0$, $t>0$, and $\rho=\rho(t)\in (0,\rho_0]$. By Plancherel's theorem and the monotonicity property of $s_\mu$ w.r.t.\ $\mu$ we have
\begin{align}
(2\pi)^{d}|u(t,\cdot)|_2^2 & = |\tilde{u}(t,\cdot)|_2^2= \int_{\iR^d} \tilde{Z}(t,\xi)^2 |\widetilde{u_0}(\xi)|^2\,d\xi
\ge \int_{B_\rho} \tilde{Z}(t,\xi)^2 |\widetilde{u_0}(\xi)|^2\,d\xi\nonumber\\
& = \int_{B_\rho} s(t,|\xi|^2)^2 |\widetilde{u_0}(\xi)|^2\,d\xi \ge s(t,\rho^2)^2 \int_{B_\rho}|\widetilde{u_0}(\xi)|^2\,d\xi\nonumber\\
& = s(t,\rho^2)^2 \rho^d \left(\rho^{-d}\int_{B_\rho}|\widetilde{u_0}(\xi)|^2\,d\xi \right).
\label{basiclowerest}
\end{align}
 Recall that $u_0\in L_1(\iR^d)\cap L_2(\iR^d)$ implies that $\widetilde{u_0}\in C_0(\iR^d)\cap L_2(\iR^d)$.

 By the assumption
 $\widetilde{u_0}(0)\neq 0$, we may choose $\rho_0$ so small that we get an estimate
 \begin{equation} \label{c1est}
 \rho^{-d}\int_{B_\rho}|\widetilde{u_0}(\xi)|^2\,d\xi\ge c_1 \quad \mbox{for all}\;\rho\in (0,\rho_0],
 \end{equation}
 with some constant $c_1>0$. Using (\ref{c1est}) and the lower estimate in (\ref{smubounds}) we deduce from (\ref{basiclowerest}) that
 \begin{equation} \label{basiclowerest2}
(2\pi)^{d} |u(t,\cdot)|_2^2 \ge \,\frac{c_1 \rho^d}{\big(1+\rho^2 k(t)^{-1}\big)^2}.
 \end{equation}

We first choose
 \[
 \rho=\rho(t)=\,\frac{\rho_0}{\big(1+k(t)^{-1}\big)^{1/2}}.
 \]
With this choice we have $\rho(t)^2 k(t)^{-1}\le \rho_0^2$ and thus
\begin{align*}
(2\pi)^{d}|u(t,\cdot)|_2^2 & \ge \,\frac{c_1 \rho^d}{\big(1+\rho_0^2\big)^2}
=\,\frac{c_1 \rho_0^d}{\big(1+\rho_0^2\big)^2}\,\left(\frac{k(t)}{1+k(t)}\right)^{
\frac{d}{2}}\\
& \ge \,\frac{c_1 \rho_0^d \, k(t)^{d/2}}{\big(1+\rho_0^2\big)^2\,\big(1+k(t_0)\big)^{d/2}},\quad t\ge t_0>0,
\end{align*}
since $k$ is nonincreasing. This shows that in this situation we have
$|u(t)|_2 \gtrsim k(t)^{d/4}$ for $t\ge 1$.

Let us next choose $\rho=\rho(t)=\rho_0$. Then (\ref{basiclowerest2}) yields
\[
(2\pi)^{d}|u(t,\cdot)|_2^2 \ge \,\frac{c_1 \rho_0^d}{\big(1+\rho_0^2 k(t)^{-1}\big)^2}\,
=\,\frac{c_1 \rho_0^d\,k(t)^2}{\big(k(t_0)+\rho_0^2\big)^2},\quad t\ge t_0>0,
\]
and hence $|u(t)|_2 \gtrsim k(t)$ for $t\ge 1$. \hfill $\square$

${}$

We turn now to upper estimates.
\begin{satz} \label{thmupper}
Let $d\in \iN\setminus \{4\}$, $u_0\in L_1(\iR^d)\cap L_2(\iR^d)$, and assume that $(k,l)\in \mathcal{PC}$.
Suppose that $u$ is given by (\ref{gensol}). Then there holds
\[
|u(t)|_2 \lesssim \,\big[(1\ast l)(t)\big]^{-\min\{1,\frac{d}{4}\}},\quad t> 0.
\]
\end{satz}
{\em Proof.} Suppose that $d\le 3$. By Plancherel's theorem and the upper estimate in (\ref{smubounds}) we have
\begin{align}
(2\pi)^{d}|u(t,\cdot)|_2^2 & = |\tilde{u}(t,\cdot)|_2^2= \int_{\iR^d} s(t,|\xi|^2)^2 |\widetilde{u_0}(\xi)|^2\,d\xi
\le |\widetilde{u_0}|_\infty^2 \int_{\iR^d} s(t,|\xi|^2)^2 \,d\xi \nonumber\\
& \le \,C|u_0|_1^2 \int_{\iR^d}
\,\frac{d\xi}{\big(1+|\xi|^2 \,(1\ast l)(t)\big)^2}=\,C_1|u_0|_1^2 \int_0^\infty
\,\frac{r^{d-1}\,dr}{\big(1+r^2 \,(1\ast l)(t)\big)^2}\nonumber\\
& = \,\frac{C_1|u_0|_1^2}{\big((1\ast l)(t)\big)^{d/2}} \int_0^\infty
\,\frac{\rho^{d-1}\,d\rho}{\big(1+\rho^2\big)^2}=\,\frac{C_2|u_0|_1^2}{\big((1\ast l)(t)\big)^{d/2}}, \label{upper1}
\end{align}
which shows that $|u(t)|_2 \lesssim (1\ast l)(t)\big)^{-d/4}$, $t>0$. Observe that this estimate breaks down for $d\ge 4$, since the last integral becomes infinite.

Suppose now that $d\ge 5$. By interpolation it follows from $u_0\in L_1(\iR^d)\cap L_2(\iR^d)$ that $u_0\in L_{\frac{2d}{d+4}}(\iR^d)$. Since
\[
\frac{1}{\frac{2d}{d+4}}\,-\,\frac{1}{2}\,=\,\frac{2}{d} \quad \mbox{and}\;\frac{2d}{d+4}>1,
\]
the Hardy-Littlewood-Sobolev theorem on fractional integration, see e.g.\ \cite[Thm.\ 6.1.3]{Graf}, implies
$(-\Delta)^{-1} u_0\in L_2(\iR^d)$. Using this property, Plancherel's theorem, and the upper estimate in (\ref{smubounds}) we may estimate as follows.
\begin{align*}
(2\pi)^{d}|u(t,\cdot)|_2^2 & = \int_{\iR^d} |\xi|^4 s(t,|\xi|^2)^2 \big||\xi|^{-2}\widetilde{u_0}(\xi)\big|^2\,d\xi\\
& \le \,\frac{1}{\big((1\ast l)(t)\big)^2}\,\int_{\iR^d} \frac{|\xi|^4\big((1\ast l)(t)\big)^2}{\big(1+|\xi|^2 \,(1\ast l)(t)\big)^2}\,\big||\xi|^{-2}\widetilde{u_0}(\xi)\big|^2\,d\xi\\
& \le \,\frac{1}{\big((1\ast l)(t)\big)^2}\,\int_{\iR^d} \big||\xi|^{-2}\widetilde{u_0}(\xi)\big|^2\,d\xi
=\,\frac{(2\pi)^{d}}{\big((1\ast l)(t)\big)^2}\,|(-\Delta)^{-1} u_0|_2^2.
\end{align*}
Hence $|u(t)|_2 \lesssim (1\ast l)(t)\big)^{-1}$, $t>0$. The theorem is proved. \hfill $\square$
\begin{bemerk1} \label{importantremark}
{\em Comparing the results from Theorem \ref{theoremlower} and Theorem \ref{thmupper} we see that
the estimate from Theorem \ref{thmupper} is in general optimal provided that
\begin{equation} \label{optimalitycond}
(1\ast l)(t)^{-1}\le
C k(t),\quad  t\ge T,
\end{equation}
where $C,T$ are some positive constants. This property of the kernels $k,l$ is satisfied in many
important examples, in particular in the time-fractional and ultraslow diffusion case, see Examples \ref{beifrac} and
\ref{beislow} below. However, there exist examples where (\ref{optimalitycond}) is violated, see e.g.\ Example \ref{exotic}, and thus Theorem \ref{thmupper} may not give the optimal decay rate. Recall that
a key ingredient in the proof of Theorem \ref{thmupper} is the upper bound
\begin{equation} \label{smuup}
s_\mu(t)\le \,\frac{1}{1+\mu \,(1\ast l)(t)}\,,\quad
\;t\ge 0,\\\mu\ge 0.
\end{equation}
If (\ref{optimalitycond}) is violated, then (\ref{smuup}) may not be sharp w.r.t.\ the decay rate. In fact, this is the case
in Example \ref{exotic} below. To improve the decay estimate from Theorem \ref{thmupper} one has to upgrade the estimate for $s_\mu$.

Suppose we have an estimate
\begin{equation} \label{smubetter}
s_\mu(t)\le \,\frac{C}{1+\mu \,\psi(t)}\,,\quad
\;t\ge 0,\,\mu\ge 0,
\end{equation}
where $\psi:\iR_+\rightarrow \iR_+$ is a continuous function (typically nondecreasing) and the constant $C$
is independent of $\mu$. Then it follows from the proof of Theorem \ref{thmupper} that
\begin{equation} \label{improvedestimate}
|u(t)|_2 \lesssim \,\big[\psi(t)\big]^{-\min\{1,\frac{d}{4}\}},\quad t> 0.
\end{equation}
}
\end{bemerk1}

\begin{bei} \label{beifrac}
The time-fractional case. {\em We consider the pair
\[
(k,l)=(g_{1-\alpha},g_\alpha),\quad \mbox{where}\;\alpha\in (0,1).
\]
Recall that the Laplace transform of $g_\beta$, $\beta>0$, is given by $\widehat{g_\beta}(z)=z^{-\beta}$, Re$\,z>0$, and so it is easy to see that $g_{\beta_1}\ast
g_{\beta_2}=g_{\beta_1+\beta_2}$ for all $\beta_1,\beta_2>0$. In particular
$(k,l)\in \mathcal{PC}$. We further have
\[
(1\ast l)(t)=(1\ast g_\alpha)(t)=g_{1+\alpha}(t)=\frac{t^\alpha}{\Gamma(1+\alpha)}
\]
and so both $k(t)$ and $(1\ast l)(t)^{-1}$ are of the form $c\,t^{-\alpha}$ with some constant $c>0$.
Invoking Theorem \ref{thmupper} reproduces the estimates from Corollary \ref{L2estimate} for noncritical $d$, that is,
\[
|u(t)|_2 \lesssim t^{-\alpha\min\{1,\frac{d}{4}\}},\quad t> 0,\;d\in \iN\setminus\{4\}.
\]
Theorem \ref{theoremlower} shows that this estimate is optimal in the general case.
}
\end{bei}
\begin{bei}
A sum of two fractional derivatives. {\em Let $0<\alpha<\beta<1$ and
\[
k(t)=g_{1-\alpha}(t)+g_{1-\beta}(t), \quad t>0.
\]
Clearly, $k$ is completely monotone and $k(0+)=\infty$, and so by Theorem 5.4 in Chapter 5 of \cite{GLS}, the kernel
$k$ has a resolvent $l\in L_{1,loc}(\iR_+)$ of the first kind,
that is $k\ast l=1$ on $(0,\infty)$, and this resolvent is
completely monotone as well. In particular $(k,l)\in {\cal PC}$. Observe that
\[
\widehat{1\ast l}\,(z)=\,\frac{1}{z}\,\frac{1}{z^\alpha+z^\beta}\,\sim \frac{1}{z^{1+\alpha}}\quad \mbox{as}\;z\to 0,
\]
which yields $(1\ast l)(t)\sim g_{1+\alpha}(t)$ as $t\to \infty$, by the Karamata-Feller Tauberian theorem, see \cite{Feller}.
It follows that
there is $T_1>0$ such that $(1\ast l)(t)\ge \frac{1}{2}g_{1+\alpha}(t)$ for all $t\ge T_1$.
Applying Theorem \ref{thmupper} then yields
\[
|u(t)|_2 \lesssim t^{-\alpha\min\{1,\frac{d}{4}\}},\quad t\ge T_1,\;d\in \iN\setminus\{4\}.
\]
Theorem \ref{theoremlower} shows that this estimate is optimal. We see that the decay rate is determined
by the fractional derivative of lower order. These considerations extend trivally to kernels $k(t)=\sum_{j=1}^m \delta _j g_{1-\alpha_j}(t)$ with $\delta_j>0$ and
$0<\alpha_1<\alpha_2<\ldots<\alpha_m<1$.
}
\end{bei}
\begin{bei} \label{beislow}
The ultraslow diffusion case. {\em We consider the pair
\[
k(t)=\int_0^1 g_\beta(t)\,d\beta,\quad l(t)=\int_0^\infty \,\frac{e^{-st}}{1+s}\,{ds},\quad t>0.
\]
Both kernels are nonnegative and nonincreasing, and there holds (see \cite[Example 6.5]{VZ})
\[
\hat{k}(z)=
\,\frac{z-1}{z\log z},\quad  \hat{l}(z)=\,\frac{\log z}{z-1}\,,\quad  \mbox{Re}\,z>0.
\]
Thus $(k,l)\in \mathcal{PC}$. There exists a number $T_1>1$ such that
\[
\frac{1}{2k(t)}\,\le \log t\le 2(1\ast l)(t),\quad t\ge T_1,
\]
see \cite[Example 6.5]{VZ}. This together with Theorem \ref{thmupper} yields the logarithmic decay estimate
\[
|u(t)|_2 \lesssim (\log t)^{-\min\{1,\frac{d}{4}\}},\quad t\ge T_1,
\]
which is optimal, by Theorem \ref{theoremlower}.
}
\end{bei}
\begin{bei} \label{exotic}
Switching the kernels from the previous example. {\em We consider now the pair
\[
k(t)=\int_0^\infty \,\frac{e^{-st}}{1+s}\,{ds},\quad l(t)=\int_0^1 g_\beta(t)\,d\beta,\quad ,\quad t>0.
\]
From the previous considerations we know already that $(k,l)\in \mathcal{PC}$.
The kernel $k(t)$ in this example behaves like $t^{-1}$ as $t\to \infty$,
whereas $(1\ast l)(t)\sim t/\log t$ as $t\to \infty$, see \cite[Example 6.6]{VZ}. Thus $k(t)$ decays faster than $(1\ast l)(t)^{-1}$, so that there is
a gap between the decay rates provided by
Theorem \ref{theoremlower} and Theorem \ref{thmupper}. We claim that
\[
|u(t)|_2 \lesssim t^{-\min\{1,\frac{d}{4}\}},\quad t> 0,\;d\in \iN\setminus\{4\},
\]
which is optimal by Theorem \ref{theoremlower}. What is interesting here is that for $d\le 3$ the decay rate is the same as for the heat equation!

The claim follows from the estimate
\begin{equation} \label{smuspecial}
s_\mu(t)\le \frac{C}{1+\mu t},\quad t\ge 0,\;\mu\ge 0,
\end{equation}
which has been shown in \cite[Example 6.6]{VZ}, and from (\ref{improvedestimate}) in Remark \ref{importantremark}.
The proof of (\ref{smuspecial}) in \cite{VZ} is quite involved and makes use of Laplace transform methods.
}
\end{bei}
\section{$L_r$-estimates and the critical dimension case} \label{SecLR}
In this section we continue the study of the general subdiffusion problem (\ref{subdiffgen}), (\ref{initcond}) under the
assumption that the kernel $k$ is of type $\cal{PC}$.
\begin{lemma} \label{aux1}
Assume that $(k,l)\in \cal{PC}$. Let $t\ge 0$ be fixed. Then the function $\mu\mapsto s(t,\mu)$
belongs to $C^\infty(\iR_+)$ and
\begin{equation} \label{smuCM}
(-1)^j \partial_\mu^j s(t,\mu)\ge 0\quad \mbox{for all}\;j\in \iN_0,
\end{equation}
in particular $s(t,\mu)$ is completely monotone w.r.t.\ $\mu$.
Moreover,
\begin{equation} \label{smuHest}
\mu^j \big|\partial_\mu^j s(t,\mu)\big|\le 2^j j!\,s(t,\frac{\mu}{2})\quad \mbox{for all}\;j\in \iN_0.
\end{equation}
\end{lemma}
{\em Proof.} Recall that $s(t,\mu)=s_\mu(t)$ solves the equation
\[
s_\mu(t)+\mu (s_\mu\ast l)(t)=1,\quad t,\,\mu\ge 0.
\]
Since $\mu$ merely appears as coefficient in front of the second term, it is clear that the dependence of the solution $s_\mu(t)$ on the parameter $\mu$ is $C^\infty$. Differentiating w.r.t.\ $\mu$ gives
\[
\partial_\mu s_\mu+\mu (\partial_\mu s_\mu \ast l)+s_\mu\ast l=0,
\]
which is equivalent to the integro-differential equation
\begin{equation} \label{deriv1}
\frac{d}{dt}\left(k\ast \partial_\mu s_\mu\right)(t)+\mu \partial_\mu s_\mu(t)=-s_\mu(t),\quad t>0,\quad \partial_\mu s_\mu(0)=0.
\end{equation}
The property $\partial_\mu s_\mu(0)=0$ follows from the fact that $s_\mu(0)=1$ for all $\mu\ge 0$. Note also
that $\partial_\mu s_\mu|_{\mu=0}=-(1\ast l)(t)$. Differentiating (\ref{deriv1}) w.r.t.\ $\mu$ leads to
\[
\frac{d}{dt}\left(k\ast \partial_\mu^2 s_\mu\right)(t)+\mu \partial_\mu^2 s_\mu(t)=-2\partial_\mu s_\mu(t),\quad t>0,\quad \partial_\mu^2 s_\mu(0)=0.
\]
Differentiating further, a simple induction argument shows that
\begin{equation} \label{deriv2}
\frac{d}{dt}\left(k\ast \partial_\mu^j s_\mu\right)(t)+\mu \partial_\mu^j s_\mu(t)=-j \partial_\mu^{j-1} s_\mu (t),\quad t>0,\quad \partial_\mu^j s_\mu(0)=0,\quad j\in\iN.
\end{equation}
Assertion (\ref{smuCM}) follows then by means of induction from (\ref{deriv2}) and the fact that the solution $v$ of
\[
\frac{d}{dt}\left(k\ast v\right)(t)+\mu v(t)=f(t),\quad t>0,\quad v(0)=0,
\]
is nonnegative, whenever $f\in L_{1,\,loc}(\iR_+)$ enjoys this property, see e.g.\ \cite{VZ} for the latter property.

To see (\ref{smuHest}), we apply Taylor's theorem to the function $\mu\mapsto s(t,\mu)$, thereby obtaining that for every
$n\in \iN$
\begin{equation} \label{Taylor}
s(t,\frac{\mu}{2})=\sum_{j=0}^n \frac{\partial_\mu^j s(t,\mu)}{j!}\,(-\frac{\mu}{2})^j+\frac{\partial_\mu^{n+1} s(t,\eta)}{(n+1)!}\,(-\frac{\mu}{2})^{n+1},
\end{equation}
for some $\eta\in (\frac{\mu}{2},\mu)$. In view of (\ref{smuCM}) every summand on the right-hand side of (\ref{Taylor}) is nonnegative. This implies
\[
s(t,\frac{\mu}{2})\ge \frac{\partial_\mu^j s(t,\mu)}{j!}\,(-\frac{\mu}{2})^j \quad \mbox{for all}\,j\le n,
\]
which in turn yields (\ref{smuHest}). \hfill $\square$
\begin{lemma} \label{aux2}
Let $(k,l)\in \cal{PC}$ and $t\ge 0$ be fixed. Let $\kappa\in (0,1]$ and set $\psi_\kappa(\mu)=\mu^\kappa s(t,\mu)$, $\mu> 0$.
Then $\psi_\kappa\in C^\infty((0,\infty))$ and for every $n\in \iN$ there exists a constant $C(n)>0$ such that
\[
\mu^n |\psi_\kappa^{(n)}(\mu)|\,[ (1\ast l)(t)]^\kappa \le C(n),\quad \mu> 0.
\]
\end{lemma}
{\em Proof.} By Leibniz' formula for the $n$th derivative of a product of two functions we have
\[
\mu^n\psi_\kappa^{(n)}(\mu)=\sum_{j=0}^n \binom{n}{j} \big(\mu^j\partial_\mu^j s(t,\mu)\big)\cdot \big(\mu^{n-j}\partial_\mu^{n-j}(\mu^\kappa)\big),
\]
and thus by Lemma \ref{aux1} and (\ref{smubounds}),
\begin{align*}
\mu^n |\psi_\kappa^{(n)}(\mu)| & \,\le \sum_{j=0}^n \binom{n}{j} \big(\mu^j|\partial_\mu^j s(t,\mu)|\big)\cdot \big(\mu^{n-j}|\partial_\mu^{n-j}(\mu^\kappa)|\big)\\
& \,\le \tilde{c}(n) \mu^\kappa \sum_{j=0}^n \binom{n}{j} 2^j j!\, s(t,\frac{\mu}{2})\\
& \,\le C(n)\,\frac{\mu^\kappa}{1+\mu (1\ast l)(t)}\,\le \, \,\frac{C(n)}{ [(1\ast l)(t)]^\kappa}.
\end{align*}
This proves the lemma. \hfill $\square$
\begin{lemma} \label{aux3}
Let $(k,l)\in \cal{PC}$, $\delta\in (0,1]$, and $t\ge 0$ be fixed. Let $m_0(\xi)=\psi_\delta(|\xi|^2)=|\xi|^{2\delta} s(t,|\xi|^2)$,
$\xi\in \iR^d$. Then $m_0\in C^\infty((0,\infty)^d)$ and for any multiindex $\beta=(\beta_1,\ldots,\beta_d)\in \iN_0^d$ the partial derivative $\partial_{\xi}^\beta m_0$ of order $|\beta|=\sum_{i=1}^d \beta_i$ is a sum of finitely many terms of the form
\begin{equation} \label{summandform}
c(\beta)\cdot \psi_\delta^{(j)}\left(|\xi|^2\right)\cdot \prod_{i=1}^d \xi_i^{\gamma_i}\quad
\mbox{with}\;\;\frac{|\beta|}{2}
\le j \le |\beta|\;\;\mbox{and}\;\;\sum_{i=1}^d \gamma_i=2j-|\beta|,
\end{equation}
where $c(\beta)>0$.

Moreover, the function $m(\xi):=m_0(\xi) [(1\ast l)(t)]^\delta$ satisfies
Mihlin's condition with a constant that is uniform w.r.t.\ $t\ge 0$, that is
there exists $M=M(d)>0$ such that
\[
|\xi|^{|\beta|} \big|\partial_\xi^\beta m(\xi)\big|\le M,\quad \xi\in \iR^d\setminus\{0\},\;\;|\beta|\le \left[\frac{d}{2}\right]+1.
\]
\end{lemma}
{\em Proof.} The assertion on the structure of $\partial_{\xi}^\beta m_0$ can be proved by induction over $|\beta|$.
If $|\beta|=0$, then $\partial_{\xi}^\beta m_0(\xi)=m_0(\xi)=\psi(|\xi|^2)$, which is of the desired form with $j=0$.
Suppose now that the assertion is true for all $\beta\in\iN_0^d$ of the same fixed order $b:=|\beta|\in \iN_0$. Let
$\beta'\in \iN_0^d$ with $|\beta'|=b+1$. Then $\partial_{\xi}^{\beta'} m_0=\partial_{\xi_l}\partial_{\xi}^\beta m_0$
for some $\beta\in \iN_0^d$ with $|\beta|=b$ and some $l\in \{1,\ldots,d\}$. By the induction hypothesis, $\partial_{\xi_l}\partial_{\xi}^\beta m_0$ is a finite sum of first order partial derivatives w.r.t.\ $\xi_l$ of terms of the
form described in (\ref{summandform}). Let us consider such a term. If $\gamma_l=0$ we have
\[
\partial_{\xi_l}\left[ \psi_\delta^{(j)}\left(|\xi|^2\right)\cdot \prod_{i=1}^d \xi_i^{\gamma_i}\right]=
\psi_\delta^{(j+1)}\left(|\xi|^2\right)\cdot 2\xi_l \cdot \prod_{i=1}^d \xi_i^{\gamma_i},
\]
whereas in case $\gamma_l>0$ we obtain
\begin{align}
\partial_{\xi_l}\left[ \psi_\delta^{(j)}\left(|\xi|^2\right)\cdot \prod_{i=1}^d \xi_i^{\gamma_i}\right]=
\psi_\delta^{(j+1)}\left(|\xi|^2\right)\cdot 2\xi_l \cdot \prod_{i=1}^d \xi_i^{\gamma_i}+ \psi_\delta^{(j)}\left(|\xi|^2\right)\cdot
\gamma_l \xi_l^{\gamma_l-1}\prod_{i=1,i\neq l}^d \xi_i^{\gamma_i}.
\label{inductionderiv}
\end{align}
The first term on the right-hand side of (\ref{inductionderiv}) has the desired form, since with $\gamma_i':=\gamma_i$, $
i\neq l$ and $\gamma_l'=\gamma_l+1$, we have by the induction hypothesis
\[
0\le \sum_{i=1}^d \gamma'_i=\sum_{i=1}^d \gamma_i+1=2j-|\beta|+1=2(j+1)-|\beta'|.
\]
The second term has the desired form as well, since setting $\gamma_i'=\gamma_i$, $
i\neq l$ and $\gamma_l'=\gamma_l-1$, we have now
\[
0\le \sum_{i=1}^d \gamma'_i=\sum_{i=1}^d \gamma_i-1=2j-|\beta|-1=2j-|\beta'|.
\]

The second part of Lemma \ref{aux3} follows from the first one and Lemma \ref{aux2}. In fact, for any term $T(\xi)$ of the form (\ref{summandform}) Lemma \ref{aux2} yields the estimate
\begin{align*}
|\xi|^{|\beta|}|T(\xi)| & \le c(\beta)\,|\xi|^{|\beta|}\,\big|\psi_\delta^{(j)}\left(|\xi|^2\right)\big|\, \prod_{i=1}^d |\xi_i|^{\gamma_i}\\
& \le c(\beta)\,|\xi|^{2j}\big|\psi_\delta^{(j)}\left(|\xi|^2\right)\big|\cdot |\xi|^{|\beta|-2j}\,|\xi|^{\sum_{i=1}^d \gamma_i}\\
& \le  \,\frac{C(\beta,j)}{[(1\ast l)(t)]^\delta}.
\end{align*}
It is now evident that $m(\xi)$ satisfies Mihlin's condition with a constant $M$ that merely depends on the dimension $d$. \hfill $\square$

${}$

\noindent Relying on Lemma \ref{aux3}, we are now able to prove the following generalization of
Theorem \ref{Lrestimate}.
\begin{satz} \label{LrestMihlin}
Let $(k,l)\in \cal{PC}$ and suppose that $u$ is given by $u(t)=Z(t)\star u_0$, where $u_0$ is as described below.

(i) Let $d\in \iN$, $1< p<\kappa(d)$, $1<q, r<\infty$, $1+\frac{1}{r}=\frac{1}{p}+
\frac{1}{q}$, and $u_0\in  L_q(\iR^d)$. Then
\[
|u(t)|_r \lesssim \big[(1\ast l)(t)\big]^{ -\frac{ d }{2}\,\left(1-\frac{1}{p}\right)},\quad t>0.
\]
(ii) Let $d\ge 3$, $1< q,r < \infty$, $\frac{1}{r}+\frac{2}{d}=\frac{1}{q}$, and $u_0\in L_q(\iR^d)$. Then
\[
|u(t)|_r \lesssim \big[(1\ast l)(t)\big]^{ -1},\quad t>0.
\]
(iii) Let $d\ge 3$ and $u_0\in L_1(\iR^d)$. Then
\[
|u(t)|_{\frac{d}{d-2},\infty}\lesssim \big[(1\ast l)(t)\big]^{ -1} ,\quad t>0.
\]
\end{satz}
{\em Proof.} (i) Set $\delta=\frac{ d }{2}(1-\frac{1}{p})$. If $d\le 2$ it is clear that $\delta\in (0,1)$. If $d\ge 3$
we have by assumption $1<p<\kappa(d)=\frac{d}{d-2}$, which is equivalent to $\delta\in (0,1)$. With $t>0$ being fixed we write
\begin{equation} \label{decomp2}
\tilde{u}(t,\xi)=[(1\ast l)(t)]^{-\delta}\,\big(\psi_\delta(|\xi|^2) [(1\ast l)(t)]^\delta\big)\,\big(|\xi|^{-2\delta}\widetilde{u_0}(\xi)\big).
\end{equation}
By the Hardy-Littlewood-Sobolev theorem on fractional integration, see e.g.\ \cite[Thm.\ 6.1.3]{Graf},
$(-\Delta)^{-\delta}u_0\in L_r(\iR^d)$ and $|(-\Delta)^{-\delta}u_0|_r\le C(d,\delta,q)|u_0|_q$; in fact, the choice of $\delta$ and the assumption  $1+\frac{1}{r}=\frac{1}{p}+
\frac{1}{q}$ imply that
\[
\frac{1}{q}\,-\,\frac{1}{r}\,=\,\frac{2\delta}{d}\quad \mbox{and}\;\;2\delta<d.
\]
Thanks to Lemma \ref{aux3} we know that $m(\xi)=\psi_\delta(|\xi|^2) [(1\ast l)(t)]^\delta$
satisfies Mihlin's condition with a dimensional constant that is independent of $t>0$. Thus we may apply
Mihlin's multiplier theorem, see \cite[Theorem 5.2.7]{Graf}, thereby obtaining that
\[
|u(t)|_r\le C(d,r)[(1\ast l)(t)]^{-\delta}|(-\Delta)^{-\delta}u_0|_r \lesssim [(1\ast l)(t)]^{-\delta}.
\]
This proves (i).

(ii) We consider again the decomposition (\ref{decomp2}), now setting $\delta=1$. As before we see that the Hardy-Littlewood-Sobolev theorem implies $(-\Delta)^{-1}u_0\in L_r(\iR^d)$. The assertion follows then from
Lemma \ref{aux3} with $\delta=1$ and Mihlin's multiplier theorem.

(iii) We know already that $m(\xi)=\psi_1(|\xi|^2) [(1\ast l)(t)]$ is an $L_r(\iR^d)$-Fourier multiplier for all $r\in (1,\infty)$ with a constant that only depends on $r$ and $d$, that is, the operator $T$ defined by $Tf=\mathcal{F}^{-1}(m \mathcal{F}f)$
($\mathcal{F}$ denoting the Fourier transform) on a suitable dense subset of $L_r(\iR^d)$ is $L_r(\iR^d)$-bounded, thus extends to an operator $T\in\mathcal{B}(L_r(\iR^d))$, and $|T|_{\mathcal{B}(L_r)}\le M(d,r)$. The weak $L_r$-spaces
can be obtained from the strong ones by real interpolation. Assuming $1<r<\infty$ we may choose $r_1\in (1,r)$,
$r_2\in (r,\infty)$, and $\theta\in (0,1)$ such that $\frac{1}{r}=\frac{1-\theta}{r_1}+\frac{\theta}{r_2}$. By \cite[Theorem 1.18.2]{Triebel}, we then have $(L_{r_1},L_{r_2})_{\theta,\infty}=L_{r,\infty}$. It follows that $T\in \mathcal{B}(L_{r,\infty}(\iR^d))$, with a norm bound that only depends on $r$ and $d$.

We choose $r=\frac{d}{d-2}$. Then $1-\frac{1}{r}=\frac{2}{d}$, and
the Hardy-Littlewood-Sobolev theorem (\cite[Thm.\ 6.1.3]{Graf}) implies that $(-\Delta)^{-1}u_0\in L_{r,\infty}(\iR^d)$.
Note that $u_0\in L_1(\iR^d)$ and so we only get an estimate in a weak $L_r$-space.
The assertion now follows from (\ref{decomp2}) with $\delta=1$ and the fact that $T\in \mathcal{B}(L_{r,\infty}(\iR^d))$, with a norm bound that is independent of $t>0$. \hfill $\square$

${}$

\noindent The following result generalizes Corollary \ref{Lrkorollar} and is a direct consequence of the previous theorem. The proof is analogous to the one of  Corollary \ref{Lrkorollar}.
\begin{korollar} \label{Lrkorollargeneral}
Let $(k,l)\in \cal{PC}$, $1<r<\infty$, $d\in \iN$, $u_0\in L_1(\iR^d)\cap L_r(\iR^d)$, and $u(t)=Z(t)\star u_0$. Then
\begin{align*}
|u(t)|_r \,& \lesssim\,\big[(1\ast l)(t)\big]^{ -\frac{ d }{2}\,\left(1-\frac{1}{r}\right)},\quad t>0,\;\; \mbox{if}\;\;d<3,\\
|u(t)|_r \,& \lesssim\,\big[(1\ast l)(t)\big]^{ -\frac{ d }{2}\,\left(1-\frac{1}{r}\right)},\quad t>0,\;\; \mbox{if}\;\;d\ge 3\;\mbox{and}\;d<\frac{2r}{r-1},\\
|u(t)|_{r,\infty} \,& \lesssim\,\big[(1\ast l)(t)\big]^{-1},\quad t>0,\;\; \mbox{if}\;\;d\ge 3\;\mbox{and}\;d=\frac{2r}{r-1},\\
|u(t)|_r \,& \lesssim\,\big[(1\ast l)(t)\big]^{-1},\quad t>0,\;\; \mbox{if}\;\;d\ge 3\;\mbox{and}\;d>\frac{2r}{r-1}.
\end{align*}
\end{korollar}

${}$

\noindent As a special case we obtain the expected weak $L_2$-decay estimate in the case of the critical dimension $d=4$, which was missing in Theorem \ref{thmupper}.
\begin{korollar} \label{L2estcriticaldimension}
Let $d=4$, $u_0\in L_1(\iR^d)$, and assume that $(k,l)\in \mathcal{PC}$.
Suppose that $u$ is given by $u(t)=Z(t)\star u_0$. Then
\[
|u(t)|_{2,\infty} \lesssim\,\big[(1\ast l)(t)\big]^{-1},\quad t> 0.
\]
\end{korollar}
\begin{bemerk1}
{\em Similarly as in Section \ref{SecL2} the decay results from this section might not be optimal if the condition (\ref{optimalitycond}) is violated, but can be improved provided one has an estimate (\ref{smubetter}) where $\psi(t)$ increases faster than
$(1\ast l)(t)$ as $t\to \infty$. In this case, the above statements from Lemma \ref{aux2}, Lemma \ref{aux3}, Theorem
\ref{LrestMihlin}, Corollary \ref{Lrkorollargeneral}, and Corollary \ref{L2estcriticaldimension} remain valid when
$(1\ast l)(t)$ is replaced by $\psi(t)$.
}
\end{bemerk1}
\section{Decay estimates via the energy method}
In this section we study a time-fractional diffusion equation with a more general
elliptic operator in divergence form. To be precise, we consider the problem
\begin{align}
\partial_t^\alpha (u-u_0)-\mbox{div}\,\big(A(t,x) \nabla u\big) & =0,\quad t>0,\,x\in \iR^d, \label{Afracdiff}\\
u|_{t=0} & = u_0,\quad x\in \iR^d. \label{Afracdiffini}
\end{align}
Here, we assume $\alpha\in (0,1)$, $u_0\in L_1(\iR^d) \cap L_2(\iR^d)$, and
\begin{itemize}
\item [{\bf ($\mathcal{H}$)}] $A\in L_\infty((0,T)\times B_R;\iR^{d\times
d})$ for all $T, R>0$, and $\exists \nu>0$ such that
\[
\big(A(t,x)\xi|\xi\big)\ge \nu|\xi|^2,\quad \mbox{for
a.a.}\;(t,x)\in (0,\infty)\times \iR^d,\,\mbox{and all}\,\xi\in
\iR^d,
\]
\end{itemize}
where $B_R=B_R(0)$ denotes the ball of radius $R$ centered at $0$.
Let $\oH^1_2(B_R)=\overline{C_0^\infty(B_R)}\,{}^{H^1_2(B_R)}$.
We say that a function $u:(0,\infty)\times \iR^d \rightarrow \iR$ is a global weak solution of (\ref{Afracdiff}), (\ref{Afracdiffini}) if for any $T, R>0$,
\begin{align*}
u|_{(0,T)\times B_R} &\in \{\,v\in
L_2([0,T];H^1_2(B_R))\;
\mbox{such that}\;\\
&\;\;g_{1-\alpha}\ast v\in C([0,T];L_2(B_R)),
\;\mbox{and}\;(g_{1-\alpha}\ast v)|_{t=0}=0\},
\end{align*}
and for any test function
\[
\eta\in H^1_2([0,T];L_2(B_R))\cap
L_2([0,T];\oH^1_2(B_R))
\]
with $\eta|_{t=T}=0$ there holds
\[
\int_{0}^{T} \int_{B_R} \Big(-\eta_t [g_{1-\alpha}\ast (u-u_0)]+
(A\nabla u|\nabla \eta)\Big)\,dx\,dt=0.
\]
\begin{satz} \label{energydecay}
Let $d\in \iN$, $\alpha\in (0,1)$, $u_0\in L_1(\iR^d)\cap L_{2}(\iR^d)$, and suppose that condition {\bf ($\mathcal{H}$)} is satisfied. Let $u$ be a global weak solution of (\ref{Afracdiff}), (\ref{Afracdiffini}), and assume in addition that
\begin{equation} \label{condinfinity}
|\nabla u|,\,u^2\in L_{1,loc}([0,\infty);L_1(\iR^d)).
\end{equation}
Then
\begin{equation} \label{thmEDecay1}
|u(t)|_{2} \lesssim t^{-\frac{\alpha d}{d+4}},\quad t>0.
\end{equation}
\end{satz}
An important tool in our proof is the so-called
{\em fundamental identity} for integro-differential
operators of the form $\frac{d}{dt}(k\ast \cdot)$, cf.\ also
\cite{Za1}. It can be viewed as the analogue to the chain rule
$(H(u))'=H'(u)u'$.
\begin{lemma} \label{FILemma}
Let $T>0$ and $U$ be an open subset of $\iR$. Let further $k\in
H^1_1([0,T])$, $H\in C^1(U)$, and $u\in L_1([0,T])$ with $u(t)\in U$
for a.a. $t\in (0,T)$. Suppose that the functions $H(u)$, $H'(u)u$,
and $H'(u)(\dot{k}\ast u)$ belong to $L_1([0,T])$ (which is the case
if, e.g., $u\in L_\infty([0,T])$). Then we have for a.a. $t\in
(0,T)$,
\begin{align} \label{fundidentity}
H'(u(t))&\frac{d}{dt}\,(k \ast u)(t) =\;\frac{d}{dt}\,\big(k\ast
H(u)\big)(t)+
\Big(-H(u(t))+H'(u(t))u(t)\Big)k(t) \nonumber\\
 & +\int_0^t
\Big(H(u(t-s))-H(u(t))-H'(u(t))[u(t-s)-u(t)]\Big)[-\dot{k}(s)]\,ds.
\end{align}
\end{lemma}
The lemma follows from a straightforward computation. We remark that (\ref{fundidentity}) remains valid for singular kernels $k$, like e.g.\ $k=g_{1-\alpha}$ with
$\alpha\in(0,1)$, provided that $u$ is sufficiently smooth.

The following result is new and a very useful implication of the fundamental identity.
\begin{korollar} \label{convexFI} Let $T, U, k, H$, and $u$ be as in Lemma \ref{FILemma}. Let $u_0\in \iR$, and assume in addition
that $k$ is nonnegative and nonincreasing and that $H$ is convex. Then
\begin{align} \label{convexfundidentity}
H'(u(t))&\frac{d}{dt}\,\big(k \ast [u-u_0]\big)(t) \ge \;\frac{d}{dt}\,\big(k\ast
[H(u)-H(u_0)]\big)(t),\quad \mbox{a.a.}\;t\in (0,T).
\end{align}
\end{korollar}
{\em Proof.} By the fundamental identity, convexity of $H$, and the properties of $k$, we have
\begin{align*}
H'(u(t))\frac{d}{dt}\,& \big(k \ast [u-u_0]\big)(t)=H'(u(t))\frac{d}{dt}\,\big(k \ast u\big)(t)-H'(u(t))u_0 k(t)\\
& \ge \frac{d}{dt}\,\big(k\ast
H(u)\big)(t)+
\Big(-H(u(t))+H'(u(t))[u(t)-u_0]\Big)k(t)\\
& \ge \frac{d}{dt}\,\big(k\ast
H(u)\big)(t)-H(u_0)k(t)\\
& =\,\frac{d}{dt}\,\big(k\ast
[H(u)-H(u_0)]\big)(t),\quad \mbox{a.a.}\;t\in (0,T),
\end{align*}
which shows the desired inequality.
\hfill $\square$

${}$

\noindent An important consequence of Lemma \ref{FILemma} is the so-called {\em $L_p$-norm inequality} for operators of the form $\partial_t(k\ast\cdot)$, which has been established recently
in \cite{VZ}. Specializing to our situation ($p=2$) it says the following.
\begin{lemma} \label{LemmaL2IN}
Let $T>0$ and $\Omega\subset \iR^d$ be an open set. Let $k\in H^1_{1,loc}(\iR_+)$ be nonnegative and nonincreasing. Then for any $v\in L_2((0,T)\times \Omega)$ and any $v_0\in L_2(\Omega)$ there holds
\begin{equation} \label{L2norminequ}
\int_{\Omega}v\partial_t\big(k \ast [v-v_0]\big)\,dx\ge |v(t)|_{L_2(\Omega)}
\partial_t\big(k\ast \big[|v|_{L_2(\Omega)}-|v_0|_{L_2(\Omega)}\big]\big)(t),\quad
\mbox{a.a.}\,t\in (0,T).
\end{equation}
\end{lemma}
{\em Proof.} For the
reader's convenience we give a proof for (\ref{L2norminequ}), which is also simpler
than that in the more general case considered in \cite{VZ}.

By the fundamental identity, applied twice, Fubini's theorem, and the triangle inequality for the
$L_2(\Omega)$-norm we have for a.a.\ $t\in (0,T)$
\begin{align*}
\int_{\Omega}v \partial_t(k\ast v)\,dx = &\,\int_{\Omega} \Big(\frac{1}{2}\,\partial_t
(k\ast v^2)+\frac{1}{2}\,k(t)v^2\Big)\,dx\\
&\,+\int_0^t \int_{\Omega} |v(t,x)-v(t-s,x)|^2\,dx\, [-\dot{k}(s)]\,ds\\
\ge &\,\,\frac{1}{2}\,\partial_t\big(k\ast |v|^2_{L_2(\Omega)})+\,\frac{1}{2}\,k(t)
|v(t)|^2_{L_2(\Omega)}\\
&\,+\,\frac{1}{2}\,\int_0^t \Big(|v(t)|_{L_2(\Omega)}-|v(t-s)|_{L_2(\Omega)}
\Big)^2   [-\dot{k}(s)]\,ds\\
= &\,\, |v(t)|_{L_2(\Omega)} \partial_t\big(k\ast |v|_{L_2(\Omega)}\big)(t).
\end{align*}
From this and H\"older's inequality, we infer that for a.a.\ $t\in (0,T)$
\begin{align*}
\int_{\Omega}v\partial_t\big(k \ast [v-v_0]\big)\,dx & = \int_{\Omega}v\partial_t\big(k \ast v\big)\,dx-k(t)\int_\Omega vv_0\,dx\\
& \ge |v(t)|_{2} \partial_t\big(k\ast |v|_{2}\big)(t)-k(t)|v(t)|_2 |v_0|_2\\
& = |v(t)|_{2} \partial_t\big(k\ast [|v|_{2}-|v_0|_2]\big)(t).
\end{align*}
This proves the lemma. \hfill $\square$

$\mbox{}$

\noindent We are now in position to prove Theorem \ref{energydecay}.

${}$

\noindent {\em Proof.} {\bf 1. Regularized weak formulation.} For $\mu>0$, define the kernel $h_\mu\in L_{1,loc}(\iR_+)$ via the Volterra integral equation
\[
h_\mu(t)+\mu(h_\mu\ast g_\alpha)(t)=\mu g_\alpha(t),\quad t>0,
\]
and set $g_{1-\alpha,\mu}=g_{1-\alpha}\ast h_\mu$. It is well known that $g_{1-\alpha,\mu}$ is positive and nonincreasing, it belongs to $H^1_{1,loc}([0,\infty))$, and the Yosida approximation $B_n=nB(n+B)^{-1}$, $n\in \iN$, of the operator $B:=\partial_t^\alpha$ defined in an appropriate space takes the form
$B_n=\partial_t(g_{1-\alpha,n}\ast \cdot)$, and $B_n\to B$ as $n\to \infty$.
Further, $h_\mu$ is nonnegative for all $\mu>0$, and for $1\le p<\infty$ and $f\in L_p([0,T])$, we have $h_n\ast f\to f$ in $L_p([0,T])$ as $n\to \infty$. We refer to \cite{VZ1} and \cite{Za} for more background on $B_n$.
Using these properties of the kernels $g_{1-\alpha,\mu}$ one can derive an equivalent weak formulation where the singular kernel $g_{1-\alpha}$ is replaced by
the more regular kernels $g_{1-\alpha,n}$. In fact, it follows from the above definition of weak solution that for any $R>0$ and any function $\psi\in \oH^1_2(B_R)$ there holds
\begin{align} \label{regweakform}
\int_{B_R}
\Big(\psi\partial_t [g_{1-\alpha,n}\ast(u-u_0)]+\big(h_n\ast[A\nabla u]|\nabla\psi\big)\Big)\,dx=\,0,\quad
\mbox{a.a.}\,t>0,\,n\in \iN,
\end{align}
cf.\ \cite{Za}.

{\bf 2. Positive and negative part.} Denote by $y_+$ and $y_-:=[-y]_+$ the positive
and negative part, respectively, of $y\in \iR$. Appealing to \cite[Lemma 4.1]{VZ}, it follows from Step 1 that
for any $R>0$ and any nonnegative function $\psi\in \oH^1_2(B_R)$,
\begin{align} \label{regweakformpos}
\int_{B_R}
\Big(\psi\partial_t [g_{1-\alpha,n}\ast(u_{+(-)}-[u_0]_{+(-)})]+\big(h_n\ast[A\nabla (u_{+(-)})]|\nabla\psi\big)\Big)\,dx=\,0,\quad
\mbox{a.a.}\,t>0,\,n\in \iN.
\end{align}

{\bf 3. An $L_1(\iR^d)$-bound for $u(t,\cdot)$.} Letting $R>1$ we choose
a nonnegative $\psi\in C^1_0(B_R)$ such that $\psi\equiv 1$ in $B_{R-1}$ and $\psi\le 1$ as well as $|\nabla \psi|\le 2$
in $B_R$. Assuming $t\in (0,T)$ and setting $\Lambda=|A|_{L_\infty((0,T)\times B_R)}$, where $T>0$ is arbitrarily fixed,
(\ref{regweakformpos}) implies that for a.a.\ $t\in (0,T)$
\begin{align} \label{var1}
\partial_t\Big( g_{1-\alpha,n}\ast \int_{B_R}\psi u_+\,dx\Big)\le \int_{B_R}g_{1-\alpha,n}(t)
(u_0)_+ \psi\,dx+2\Lambda\int_{B_R\setminus B_{R-1}}|\nabla (u_+)|\,dx+\rho_n(t),
\end{align}
where
\[
\rho_n(t)=\int_{B_R}\big(A\nabla(u_+)-h_n\ast[A\nabla (u_{+})]\big|\nabla\psi\big)\,dx.
\]
Note that the term in brackets on the left-hand side of (\ref{var1}) vanishes
at $t=0$. Thus convolving (\ref{var1}) with $g_\alpha$ and employing that
$g_\alpha\ast g_{1-\alpha,n}=1\ast h_n$, we obtain
\begin{align*}
h_n\ast \int_{B_R}\psi u_+\,dx\le (1\ast h_n)(t)\int_{B_R}\psi (u_0)_+\,dx
+2\Lambda\,g_\alpha\ast \int_{B_R\setminus B_{R-1}}|\nabla (u_+)|\,dx+(g_\alpha\ast \rho_n)(t),
\end{align*}
Next, sending $n\to \infty$, using the approximation property of the kernels $h_n$,
and restricting to a subsequence, if necessary, we find that
\begin{align} \label{var2}
\int_{B_{R-1}}u_+\,dx\le \int_{B_R}\psi (u_0)_+\,dx
+2\Lambda\,g_\alpha\ast \int_{\iR^d\setminus B_{R-1}}|\nabla (u_+)|\,dx,
\end{align}
for a.a.\ $t\in (0,T)$. By (\ref{condinfinity}) and Young's inequality, the second term on the right-hand
side of (\ref{var2}) tends to $0$ in $L_1([0,T])$ as $R\to \infty$. Thus, sending
$R\to \infty$ in (\ref{var2}) implies
\[
|u_+(t)|_{L_1(\iR^d)}\le |(u_0)_+|_{L_1(\iR^d)},\quad \mbox{a.a.}\,t\in (0,T).
\]
The same argument leads to the corresponding estimate for $u_-$. Since $T>0$ was arbitrary, we conclude
that $|u(t)|_1\le 2|u_0|_1$ for a.a.\ $t>0$.

{\bf 4. A fractional differential inequality for the $L_2(\iR^d)$-norm.} Let $R>1$ and $\psi$ be as in Step 3. For $t>0$ we take in (\ref{regweakform}) with $\psi$ replaced by $\eta$ the test function $\eta=u\psi^2$. This gives
\begin{align} \label{step4a}
\int_{B_R}
\Big(\psi u\partial_t [g_{1-\alpha,n}\ast(\psi u-\psi u_0)]+\big(A\nabla u|\nabla(u\psi^2)\big)\Big)\,dx=\,\rho_{n, R}(t),\quad
\mbox{a.a.}\,t>0,\,n\in \iN,
\end{align}
where
\[
\rho_{n, R}(t)=\int_{B_R}\big(A\nabla u-h_n\ast[A\nabla u]\big|\nabla\eta\big)\,dx.
\]
Assume $t\in (0,T)$ and set $\Lambda=|A|_{L_\infty((0,T)\times B_R)}$, where $T>0$ is arbitrarily fixed. By assumption ($\mathcal{H}$), we may estimate pointwise a.e.\ as follows
\begin{align*}
\big(A\nabla u|\nabla(u\psi^2)\big) &\,=\,\big(A\nabla u|\nabla u)\psi^2+2\big(A\nabla u|
u\psi\nabla \psi)\\
&\,\ge \,\nu |\nabla  u|^2\psi^2-\frac{\nu}{2}\, |\nabla u|^2\psi^2-\frac{2\Lambda^2}{\nu} u^2 |\nabla \psi|^2\\
&\,\ge \,\frac{\nu}{2}\big|\nabla(u\psi)-u\nabla \psi\big|^2-\frac{2\Lambda^2}{\nu} u^2 |\nabla \psi|^2\\
&\,\ge \,\frac{\nu}{4}\,|\nabla(u\psi)|^2-\big(\frac{\nu}{2}+\frac{2\Lambda^2}{\nu}\big) u^2 |\nabla \psi|^2.
\end{align*}

The first term in (\ref{step4a}) can be estimated from below by means of Lemma \ref{LemmaL2IN}. Together with the preceding inequality we then obtain
\begin{align*}
|(\psi u)(t)|_{L_2(B_R)} &
\partial_t\big(g_{1-\alpha,n}\ast \big[|\psi u|_{L_2(B_R)}-|\psi u_0|_{L_2(B_R)}\big]\big)(t)
+\,\frac{\nu}{4}\,\int_{B_R}|\nabla(u\psi)|^2\,dx\\
\le &\,
C(\nu,\Lambda)\int_{B_R} u^2 |\nabla \psi|^2\,dx+\rho_{n, R}(t)\\
\le &\,4C(\nu,\Lambda)\int_{B_R\setminus B_{R-1}} u^2\,dx+\rho_{n, R}(t), \quad
\mbox{a.a.}\,t\in (0,T).
\end{align*}
By Nash's inequality (cf.\ \cite{Nash}, \cite{SalCoste}) and Step 3, we have for $v=\psi u$ (with $R$ being fixed)
\begin{align*}
|v(t)|_{L_2(\iR^d)}^{2+\frac{4}{d}} &\,\le C(d)|v(t)|_{L_1(\iR^d)}^{\frac{4}{d}}
|\nabla v(t)|^2_{L_2(\iR^d)}\\
&\,\le C(d) \big(2|u_0|_{L_1(\iR^d)}\big)^{\frac{4}{d}}|\nabla v(t)|^2_{L_2(\iR^d)}.
\end{align*}
Since $\psi$ vanishes outside $B_R$, we thus obtain
\begin{align*}
|v(t)|_{L_2(\iR^d)} &
\partial_t\big(g_{1-\alpha,n}\ast \big[|v|_{L_2(\iR^d)}-|\psi u_0|_{L_2(\iR^d)}\big]\big)(t)
+\,\mu |v(t)|_{L_2(\iR^d)}^{2+\frac{4}{d}}\\
&\le
\,4C(\nu,\Lambda)\int_{\iR^d\setminus B_{R-1}} u^2\,dx+\rho_{n, R}(t), \quad
\mbox{a.a.}\,t\in (0,T),
\end{align*}
for some constant $\mu=\mu(\nu,d,|u_0|_{L_1(\iR^d)})>0$. Letting $\varepsilon>0$, we next divide by
$|v(t)|_{L_2(\iR^d)}+\varepsilon$. By Corollary \ref{convexFI}, applied to
$k=g_{1-\alpha,n}$ and the convex function $H(y)=y-\varepsilon \log(y+\varepsilon)$, $y>-\varepsilon$,
with derivative $H'(y)=1-\frac{\varepsilon}{y+\varepsilon}=\frac{y}{y+\varepsilon}$
it follows that
\begin{align*}
\partial_t & \big(g_{1-\alpha,n} \ast \big[|v|_{L_2(\iR^d)}-|\psi u_0|_{L_2(\iR^d)}\big]\big)(t)+\,\mu \frac{|v(t)|_{L_2(\iR^d)}^{2+\frac{4}{d}}}{|v(t)|_{L_2(\iR^d)}+\varepsilon}\\
& \le \,
\,\frac{4C(\nu,\Lambda)}{|v(t)|_{L_2(\iR^d)}+\varepsilon}\int_{\iR^d\setminus B_{R-1}} u^2\,dx+
\frac{\rho_{n, R}(t)}{|v(t)|_{L_2(\iR^d)}+\varepsilon}\\
& \quad\,+\varepsilon \partial_t\big(g_{1-\alpha,n}\ast \big(\log(|v|_{L_2(\iR^d)}+\varepsilon)-
\log(|\psi u_0|_{L_2(\iR^d)}+\varepsilon)\big), \quad
\mbox{a.a.}\,t\in (0,T).
\end{align*}

Next, let $\varphi\in C^1([0,T])$ be nonnegative and such that $\varphi(T)=0$. We multiply the preceding inequality
by $\varphi(t)$, integrate on $(0,T)$, and integrate by parts (first and last term). Sending then $n\to \infty$, the term involving
$\rho_{n, R}$ drops. Sending next $R\to \infty$ and recalling the dependence of $\psi$ (and thus also of $v$) on $R$,
the term involving $\int_{\iR^d\setminus B_{R-1}} u^2\,dx$ goes to zero (by assumption (\ref{condinfinity})) and we
obtain
\begin{align*}
\int_0^T &\Big(-\varphi_t\, g_{1-\alpha}\ast \big[|u|_{L_2(\iR^d)}-|u_0|_{L_2(\iR^d)}\big]
+\mu\varphi \,\frac{|u|_{L_2(\iR^d)}^{2+\frac{4}{d}}}{|u|_{L_2(\iR^d)}+\varepsilon}\Big)\,dt\\
& \le \varepsilon \int_0^T \Big(-\varphi_t\, g_{1-\alpha}\ast \big(\log(|u|_{L_2(\iR^d)}+\varepsilon)-
\log(|u_0|_{L_2(\iR^d)}+\varepsilon)\big)\Big)\,dt\\
& \le 2\varepsilon |\log \varepsilon| \int_0^T |\varphi_t|g_{2-\alpha}\,dt.
\end{align*}
Sending now $\varepsilon\to 0$ the right-hand side drops, and we conclude that the $L_2(\iR^d)$-norm of $u(t,\cdot)$
satisfies the fractional differential inequality
\begin{equation} \label{basicfracdiff}
\partial_t^\alpha \big(|u|_2-|u_0|_2\big)(t)+\mu |u(t)|_2^{1+\frac{4}{d}}\le 0,\quad t\in (0,T),
\end{equation}
in the weak sense.

{\bf 5. Comparison principle and decay estimate.} By the comparison principle for time-fractional differential
equations (see \cite[Lemma 2.6 and Remark 2.1]{VZ}), (\ref{basicfracdiff}) implies that $|u(t)|_2\le w(t)$ for a.a.\
$t\in (0,T)$, where $w$ solves the equation corresponding to (\ref{basicfracdiff}), that is
\[
\partial_t^\alpha (w-w_0)(t)+\mu w(t)^\gamma= 0,\quad t>0,\quad w(0)=w_0:= |u_0|_2,
\]
where we put $\gamma={1+\frac{4}{d}}$. It is known that for $w_0>0$ there exist constants $c_1, c_2>0$ such that
\[
\frac{c_1}{1+t^{\frac{\alpha}{\gamma}}}\,\le w(t)\le\,\frac{c_2}{1+t^{\frac{\alpha}{\gamma}}},\quad t\ge 0,
\]
see \cite[Theorem 7.1]{VZ}. Since $T>0$ was arbitrary, we conclude that
\[
|u(t)|_2 \le \,w(t)\le \,\frac{c_2}{1+t^{\frac{\alpha}{\gamma}}}\,=\,\frac{c_2}{1+t^{\frac{\alpha d}{d+4}}},\quad
\mbox{a.a.}\;t>0.
\]
This finishes the proof of Theorem \ref{energydecay}. \hfill $\square$

${}$

\noindent Note that the decay rate for the $L_2$-norm in (\ref{thmEDecay1}) is strictly less than the one we obtained in
Section \ref{fractionalheat}, see (\ref{L2decay}). However, sending $d\rightarrow \infty$ the decay rate in
(\ref{thmEDecay1}) becomes $\alpha$, which is precisely the decay rate for the $L_2$-norm in (\ref{L2decay}) for $d> 4$. This phenomenon of a smaller decay rate in the variational setting does not occur in the case $\alpha=1$. In fact, proceeding similarly as above in the latter case one obtains the differential inequality
\[
\partial_t |u(t)|_{2}+\mu|u(t)|_{2}^\gamma\le 0,\quad t>0,
\]
with the same $\gamma\,(>1)$ and $\mu$ as before. This implies
\[
|u(t)|_2 \lesssim t^{-\frac{1}{\gamma-1}} = t^{-\frac{d}{4}}.
\]
On the other hand,  the $L_p(\iR^d)$-norm of the Gaussian heat kernel $H(t,x)=(4\pi t)^{-d/2}\exp(-\frac{|x|^2}{4t})$ decays as
\[
|H(t)|_p \lesssim t^{-\frac{d}{2}(1-\frac{1}{p})},\quad t>0,\;1\le p\le \infty.
\]
Assuming $u_0\in L_1(\iR^d)$, Young's inequality implies that the $L_2(\iR^d)$-norm of $u(t)=H(t)\star u_0$ decays as
\[
|u(t)|_2 \lesssim t^{-\frac{d}{2}(1-\frac{1}{2})}= t^{-\frac{d}{4}},\quad t>0,
\]
which is precisely the decay rate we obtain by means of energy estimates. It is an interesting open problem whether the decay rate
in (\ref{thmEDecay1}) is optimal in general. If this was not the case how could it be upgraded?

$\mbox{}$

\noindent {\footnotesize {\bf Jukka Kemppainen}, Mathematics Division, Dept.\ of Electr.\ and Information Engineering
Faculty of Technology, University of Oulu,
PL 4500, 90014 Oulu, Finland,
e-mail: jukemppa@paju.oulu.fi

$\mbox{}$

\noindent{\bf Juhana Siljander}, Department of Mathematics and Statistics, University of Helsinki, P.O. Box 68,
00014 Helsinki, Finland, e-mail: juhana.siljander@helsinki.fi

$\mbox{}$

\noindent {\bf Vicente Vergara}, Universidad de
Tarapac\'{a}, Instituto de Alta Investigaci\'{o}n, Antofagasta N.
1520, Arica, Chile, e-mail: vvergaraa@uta.cl

$\mbox{}$

\noindent {\bf Rico Zacher}, Martin-Luther-Universit\"at
Halle-Wittenberg, Institut f\"ur Mathematik, Theodor-Lieser-Strasse
5, 06120 Halle, Germany, e-mail: rico.zacher@mathematik.uni-halle.de

}


\begin{thebibliography}{99}
{\footnotesize
\bibitem{BS} Bjorland, C.; Schonbek, M. E.: Poincar\'e's inequality and diffusive evolution equations. Adv.\ Differential Equations {\bf 14} (2009), 241-–260.
\bibitem{BrdP} Br\"andle, C.; de Pablo, A.: Decay estimates for linear and nonlinear nonlocal heat equations. Available online at  http://arxiv.org/pdf/1312.4661v1.pdf	
\bibitem{BG}
Bouchaud, J.-Ph.; Georges, A.: Anomalous diffusion in disordered media: statistical mechanisms, models and physical applications. Phys. Rep. {\bf 195} (1990), 127--293.
\bibitem{CaffVazq11} Caffarelli, L.; Vazquez, J.\ L.: Asymptotic behaviour of a porous medium equation with fractional diffusion. Discrete Contin. Dyn. Syst. {\bf 29} (2011), 1393--1404.
\bibitem{CapuFlow}
Caputo, M.: Diffusion of fluids in porous media with memory.
Geothermics {\bf 28} (1999), 113--130.
\bibitem{CCR} Chasseigne, E.; Chaves, M.; Rossi, J.\ D.: Asymptotic behavior for nonlocal diffusion
equations. J. Math. Pures Appl. {\bf 86} (2006), 271-–291.
\bibitem{CNa} Cl\'{e}ment, Ph.; Nohel, J.A.:
Abstract linear and nonlinear Volterra equations preserving positivity.
SIAM J.\ Math.\ Anal.\ {\bf 10} (1979), 365--388.
\bibitem{CN} Cl\'{e}ment, Ph.; Nohel, J.A.: Asymptotic behavior of
solutions of nonlinear Volterra equations with completely positive
kernels. SIAM J. Math. Anal. {\bf 12} (1981), 514--534.
\bibitem{DrKl} Dr\"ager, J.; Klafter, J.:
Strong anomaly in diffusion generated by iterated maps.
Phys.\ Rev.\ Lett.\ {\bf 84} (2000), 5998--6001.
\bibitem{DZua} Duoandikoetxea, J.; Zuazua, E.: Moments, masses de Dirac et d\'ecomposition de fonctions. C.\ R.\ Acad.\ Sci.\ Paris S\'er. I Math. {\bf 315} (1992), 693--698.
\bibitem{Koch} Eidelman, S. E.; Kochubei, A. N.: Cauchy problem for
fractional diffusion equations. J. Differ. Eq. {\bf 199} (2004),
211--255.
\bibitem{Feller}
Feller, W.: {\em An introduction to probability theory and its applications.} Vol. II. Second edition. John Wiley \& Sons, Inc., New York, 1971.
\bibitem{Graf} Grafakos, L.: {\em Classical and modern Fourier analysis}. Pearson Education, 2004.
\bibitem{Grip1} Gripenberg, G.: Volterra integro-differential
equations with accretive nonlinearity. J.\ Differ.\ Eq.\ {\bf 60}
(1985), 57--79.
\bibitem{GLS} Gripenberg, G.; Londen, S.-O.; Staffans, O.:
{\em Volterra integral and functional equations.} Encyclopedia of
Mathematics and its Applications, {\bf 34}. Cambridge University
Press, Cambridge, 1990.
\bibitem{IR} Ignat, L.\ I.; Rossi, J.\ D.: Decay estimates for nonlocal problems via energy methods. J.
Math. Pures Appl. {\bf 92} (2009), 163-–187.
\bibitem{JakuDiss}
Jakubowski, V.\ G.: {\em Nonlinear elliptic-parabolic
integro-differential equations with $L_1$-data: existence,
uniqueness, asymptotics}. Dissertation, University of Essen, 2001.
\bibitem{KiSa} Kilbas, A.\ A.; Saigo, M.: {\em H-transforms: Theory and Application}. CRC Press, LLC, 2004.
\bibitem{KST} Kilbas, A.\ A.; Srivastava, H.\ M.; Trujillo, J.\ J.:
{\em Theory and applications of fractional differential
equations}. Elsevier, 2006.
\bibitem{Koch08} Kochubei, A.\ N.: Distributed order calculus and equations of
ultraslow diffusion. J. Math. Anal. Appl. {\bf 340} (2008), 252--281.
\bibitem{Koch90} Kochubei, A.\ N.: Fractional-order diffusion, Differential Equations 26 (1990) 485–492.
\bibitem{Koch11} Kochubei, A. N.: General fractional calculus, evolution equations, and renewal processes.
Integr.\ Equ.\ Oper.\ Theory {\bf 71} (2011), 583--600.
\bibitem{MZL} Ma, Y.; Zhang, F.; Changpin, L.: The asymptotics of the solutions to the anomalous
diffusion equations. Comput.\ Math.\ Appl. (2013).
\bibitem{MNV} Meerschaert, M.M.; Nane, E.; Vellaisamy, P.: Fractional Cauchy
problems on bounded domains. Ann.\ Probab.\ {\bf 37} (2009), 979--1007.
\bibitem{Metz} Metzler, R.; Klafter, J.: The random walk's guide to
anomalous diffusion: a fractional dynamics approach. Phys. Rep. {\bf
339} (2000), 1--77.
\bibitem{Metz2} Metzler, R.; Klafter, J.:
The restaurant at the end of the random walk: recent developments in the description of anomalous transport by fractional dynamics. J. Phys. A {\bf 37} (2004), R161--R208.
\bibitem{NSY} Nakagawa, J.; Sakamoto, K.; Yamamoto, M.:
Overview to mathematical analysis for fractional diffusion equations –
new mathematical aspects motivated by industrial collaboration.
Journal of Math-for-Industry {\bf 2} (2010A-10), 99--108.
\bibitem{Nash} Nash, J.: Continuity of solutions of parabolic and elliptic
equations. Amer. J. Math. {\bf 80} (1958), 931--954.
\bibitem{JanI} Pr\"uss, J.: {\em Evolutionary Integral Equations and
Applications}. Monographs in Mathematics {\bf 87}, Birkh\"auser,
Basel, 1993.
\bibitem{QS}
Quittner, P.; Souplet, Ph.: {\em Superlinear parabolic problems. Blow-up, global existence and steady states}. Birkh\"auser Verlag, Basel, 2007.
\bibitem{SalCoste} Saloff-Coste, L.: {\em Aspects of Sobolev-type
inequalities}. London Mathematical Society Lecture Note Series
{\bf 289}, Cambridge University Press, 2002.
\bibitem{SchSokBl} Schiessel, H.; Sokolov, I.\ M.; Blumen, A.:
Dynamics of a polyampholyte hooked around an obstacle. Phys.\ Rev.\ E {\bf 56} (1997), R2390--R2393.
\bibitem{SSV10} Schilling, R.; Song, R.; Vondracek, Z.: {\em Bernstein functions. Theory and applications}. Studies
in Mathematics {\bf 37}, De Gruyter, Berlin, 2010.
\bibitem{SchnWyss}
Schneider, W.\ R.; Wyss, W.: Fractional diffusion and wave
equations. J. Math. Phys.  {\bf 30} (1989), 134--144.
\bibitem{Sinai} Sinai, Y.\ G.:
The limiting behavior of a one-dimensional random walk in a random medium. Theory Probab.\ Appl.\ {\bf 27} (1982), 256--268.
\bibitem{Triebel} Triebel, H.: {\em Interpolation Theory, Function spaces, Differential operators} - 2nd, rev.\ and
enl.\ edition. Johann
Ambrosius Barth Verlag, Heidelberg, 1995.
\bibitem{Uch} Uchaikin, V.\ V.: {\em Fractional derivatives for physicists and engineers. Volume I Background and Theory.}
Nonlinear Physical Science, Springer, Heidelberg, 2013.
\bibitem{Vazq14} Vazquez, J.\ L.: Barenblatt solutions and asymptotic behaviour for a nonlinear fractional heat equation of porous medium type. To appear in J. Eur. Math. Soc, 2014.
\bibitem{VZ1} Vergara, V.; Zacher, R.: Lyapunov functions and
convergence to steady state for differential equations of fractional
order. Math. Z. {\bf 259} (2008), 287--309.
\bibitem{VZ} Vergara, V.; Zacher, R.: Optimal decay estimates for time-fractional and other non-local subdiffusion equations via energy methods. Submitted. Available online at http://arxiv.org/pdf/1310.0209.pdf
\bibitem{Za1} Zacher, R.: A De Giorgi-Nash type theorem for time fractional
diffusion equations.  Math.\ Ann.\ {\bf 356} (2013), 99--146.
\bibitem{Za} Zacher, R.: Boundedness of weak solutions to evolutionary partial
integro-differential equations with discontinuous coefficients. J.
Math. Anal. Appl. {\bf 348} (2008), 137--149.
\bibitem{ZEQ} Zacher, R.: Maximal regularity of type $L_p$ for
abstract parabolic Volterra equations. J. Evol. Equ. {\bf 5} (2005),
79--103.
\bibitem{Zua} Zuazua, E.: Large time asymptotics for heat and dissipative wave equations. Manuscript available at http://www.uam.es/enrique.zuazua, 2003.
}
\end{thebibliography}
\end{document}